\documentclass[aop,citesort,MSNbibl,seceqn,dvips]{arximspdf}

\doi{10.1214/09-AOP487}
\volume{38}
\issue{3}
\pubyear{2010}
\firstpage{1062}
\lastpage{1105}

\makeatletter
\DeclareMathAlphabet\mathcaligr{OMS}{cmsy}{m}{n}

\newtheorem{theorem}{Theorem}[section]
\newproclaim{remark}[theorem]{Remark}
\newtheorem{lemma}[theorem]{Lemma}
\newtheorem{cor}[theorem]{Corollary}
\newproclaim{defn}[theorem]{Definition}
\newtheorem{ass}{Assumption}
\newtheorem{prop}[theorem]{Proposition}

\def\shortsig{\nstop^\varepsilon_1}
\def\lcconst{\tilde{K}}
\def\stopzero{T_{\mathcaligr V}}
\def\nstop{\alpha}
\def\newstop{\beta}
\def\newtime{\zeta}
\def\newzk{Z_{(k)}}
\def\newyk{Y_{(k)}}
\def\newlk{L_{(k)}}
\def\newxk{X_{(k)}}

\makeatother

\begin{document}
\begin{frontmatter}

\title{A Dirichlet process characterization
of a~class~of~reflected diffusions}
\runtitle{Reflected diffusions and Dirichlet processes}

\begin{aug}
\author[A]{\fnms{Weining} \snm{Kang}\ead[label=e1]{weikang@andrew.cmu.edu}}
\and
\author[A]{\fnms{Kavita} \snm{Ramanan}\ead[label=e2]{kramanan@math.cmu.edu}\thanksref{t1}\corref{}}
\thankstext[1]{t1}{Supported in part by NSF Grants DMS-0406191,
DMS-0405343, CMMI-0728064.}
\runauthor{W. Kang and K. Ramanan}
\affiliation{Carnegie Mellon University and Carnegie Mellon University}
\address[A]{
Department of Mathematical\\
\quad Sciences\\
Carnegie Mellon University\\
Pittsburgh, Pennsylvania 15213\\
USA\\
\printead{e1}\\
\phantom{\textsc{E-mail: }}\printead*{e2}} 

\end{aug}

\received{\smonth{6} \syear{2008}}

%
\begin{abstract}
For a class of
stochastic differential equations with reflection for which
a certain ${\mathbb L}^p$ continuity condition holds with $p > 1$,
it is shown that any weak solution that is a strong Markov process
can be decomposed into the sum
of a local martingale and a continuous, adapted process of
zero $p$-variation.
When $p = 2$, this implies that the
reflected diffusion is a Dirichlet process.
Two examples are provided to motivate such a characterization.
The first example is a class of multidimensional reflected diffusions
in polyhedral conical domains that arise
as approximations of certain stochastic networks,
and the second example is a family of two-dimensional
reflected diffusions in
curved domains. In both cases, the reflected diffusions are shown to
be Dirichlet processes, but not semimartingales.
\end{abstract}

%

\end{frontmatter}

\section{Introduction}
\label{sec-intro}

\subsection{Background and motivation}
\label{subs-back}

This work identifies fairly general sufficient conditions under
which a reflected diffusion can be decomposed as the sum of a continuous
local martingale and a continuous adapted process of zero
$p$-variation, for some $p$ greater than one.
As motivation for such a characterization, two examples of
classes of reflected diffusions are considered.
The first example consists of
a large class of multidimensional, obliquely
reflected diffusions in polyhedral domains that
arise in applications.
Reflected diffusions in this class are shown not
to be semimartingales, but to belong to the class of
so-called Dirichlet processes.
Dirichlet processes
are processes that can be expressed (uniquely) as the sum
of a local martingale
and a continuous process that has zero quadratic variation, and
thus correspond to the case when $p = 2$ in the decomposition
mentioned above.
The second example consists of a class of two-dimensional
reflected diffusions in curved ``valley-shaped'' domains that
were first considered by Burdzy and Toby in \cite{burtoby}.
Once again, the reflected diffusion is shown
to admit a decomposition of the type mentioned above,
but in this case the magnitude of $p$ depends, in a
sense made precise in the sequel,
on the curvature of the domain.

Processes that admit a decomposition of the type mentioned above are clearly
an extension of the class of continuous semimartingales.
As is well known, semimartingales form an important class of processes
for stochastic integration,
they are stable under ${\mathcaligr C}^2$ transformations
and admit an It\^{o} change-of-variable formula.
However, there are many natural operations that lead out of the class
of semimartingales
and motivate the consideration of Dirichlet processes.
For example, ${\mathcaligr C}^1$ functionals of Brownian motion,
certain functionals
of symmetric Markov processes associated with Dirichlet forms \cite{fukbook94},
and Lipschitz functionals of a broad class of semimartingale reflected
diffusions in bounded domains \cite{roz03,rozslo00}, are all Dirichlet
processes that are in general not semimartingales.
Moreover, Dirichlet processes
exhibit many nice properties analogous to semimartingales. They admit a natural,
Doob--Meyer-type decomposition \cite{coqjakmemslo06}, they are stable
under ${\mathcaligr C}^1$ transformations (see Proposition 11 of \cite
{rusval07} and also \cite{ber89})
and there are extensions of stochastic calculus and It\^{o}'s formula
that apply to Dirichlet processes (see \cite{fol81a,folproshi95} and
Chapter 4 of \cite{rusval07}) or, more generally,
to processes that admit a decomposition as the sum of a local
martingale and
a continuous, adapted process of bounded $p$-variation, for $p \in(1,2)$
\cite{ber89}. Furthermore, the theory of rough paths
(see, e.g., \cite{FrizVictBook09} or \cite{LyonsStFlour04})
applies to
processes whose paths have bounded $p$-variation for an arbitrary
$p \in[1,\infty)$.

The theory of reflected diffusions is most well-developed for
semimartingale or symmetric
reflected diffusions.
In particular, the Skorokhod problem approach to the study
of reflected diffusions \cite{dupram99a,ram1,sko2} is automatically
limited to semimartingales,
while the Dirichlet form approach is best suited to analyze symmetric diffusions
(see, e.g., \cite{zqchen93,fukbook94}).
However, using the
submartingale formulation of Stroock and Varadhan \cite{strvarbc71} or
the extended
Skorokhod problem \cite{ram1}, it is possible to construct reflected diffusions
that are neither semimartingales nor symmetric processes
\cite{burkanram08,burtoby,ramrei1,ramrei2,wil85}. This leads
naturally to the question of
determining when these reflected diffusions are semimartingales and, if they
are not semimartingales, whether they belong to some other tractable
class of processes such as
Dirichlet processes.
There has been a substantial body of work that shows that, under certain
conditions on the domain and reflection directions (namely, the
completely-${\mathcaligr S}$ condition
and generalizations of it), the associated reflected diffusion is a
semimartingale \cite{ram1,wilsur95}.
In contrast, it has been a longstanding open problem (see Section
4(iii) of \cite{wilsur95})
to develop a theory for multidimensional reflected diffusions for which
this condition fails to hold
(some results in two dimensions can be found in \cite{burkanram08,burtoby,wil85}).
As shown in \cite{ramrei1,ramrei2}, such reflected diffusions arise
naturally as approximations
of a so-called generalized processor sharing model used in
telecommunication networks.
Thus, the development of such a theory is also of interest from the
perspective of
applications.

The first main result of this work (Theorem \ref{th-main1}) shows that
multidimensional reflected diffusions that belong to a slight
generalization of the family
of reflected diffusions obtained as approximations in \cite{ramrei1,ramrei2}
fail to be semimartingales.
In two dimensions and for the case of reflected Brownian motion,
this result follows from Theorem 5 of \cite{wil85}
(also see \cite{burkanram08} for an
alternative proof of this result).
However, the analysis in \cite{wil85} uses constructions in polar coordinates
that appear
not to be easily generalizable to higher dimensions.
We follow a different approach, which is independent of dimension and
which allows us to
establish the result for uniformly elliptic
reflected diffusions,
with possibly state-dependent diffusion coefficients, rather than just
reflected Brownian motion.

The next main result (Theorem \ref{th-main2}) shows that
reflected diffusions that belong to a broad class
admit a decomposition as the sum of a local martingale and a process of
zero $p$-variation, for some
$p > 1$.
This class consists of weak solutions to stochastic differential
equations with reflection that are Markov processes and
have locally bounded drift and dispersion
coefficients and satisfy a certain $L^p$ continuity requirement (see
Assumption \ref{as-weakcont}).
This continuity requirement is satisfied, for example, when the
associated extended Skorokhod map
is H\"{o}lder continuous,
but also holds under much weaker conditions that do not even require that
the (extended) Skorokhod map be well-defined (see Remark
\ref{rem-weakcont}).
When the extended Skorokhod map is well-defined and Lipschitz continuous,
this implies that the associated reflected diffusion is a
Dirichlet process. Using the latter result, it is shown in
Corollary \ref{cor-main2} that
the nonsemimartingale reflected diffusions considered
in Theorem \ref{th-main1} are Dirichlet processes.
Our next result concerns
the class of reflected Brownian motions introduced in \cite{burtoby},
which were shown in \cite{burkanram08} not to be semimartingales.
In Corollary \ref{cor-egburtob}, Theorem \ref{th-main2} is applied to
show that even in \mbox{cusplike} domains, the associated reflected Brownian motions are Dirichlet
processes, thus partially resolving an open question raised in \cite{burtoby}.

The paper is organized as follows. Some common notation
used throughout the paper is first summarized in Section \ref{subs-notat}.
The class of stochastic
differential equations with reflection under consideration, and the
related motivating examples, are introduced in Section \ref{sec-mainres}.
Section \ref{subs-mainres} contains a rigorous statement
of the main results;
the proof of Theorem \ref{th-main1} is presented in
Section \ref{sec-smg}, while the proofs of Theorem \ref{th-main2}
and Corollary \ref{cor-main2} are
given in Section \ref{sec-dirichlet}.
Some elementary results required in the proofs are relegated to the \hyperref[appm]{Appendix}.

\subsection{Notation}
\label{subs-notat} As usual, ${\mathbb R}_+$ or $[0,\infty)$ denote
the space
of all nonnegative reals, and
${\mathbb N}$ denotes the space of all positive integers.
Given two real numbers~$a$ and~$b$, $a \wedge b$ and $a \vee b$ denote
the minimum and maximum, respectively,
of $a$ and~$b$.
For each positive integer $J\geq1$, ${\mathbb R}^J$
denotes $J$-dimensional Euclidean space and the nonnegative orthant in
this space is denoted by ${\mathbb R}^J_+ = \{x \in{\mathbb R}^J \dvtx x_i\geq0
\mbox{ for }i = 1,\ldots, J\}.$ The Euclidean norm of $x\in
{\mathbb R}^J$ is denoted by $|x|$ and the inner product of $x,y\in
{\mathbb R}^J$ is
denoted by $\langle  x,y \rangle $.
The vectors $(e_1, e_2, \ldots, e_J)$ represent the usual orthonormal
basis for
${\mathbb R}^J$, with $e_i$ being the $i$th coordinate vector.
Given a vector $u \in{\mathbb R}^J$, $u^T$ denotes its transpose, with
analogous notation for matrices
For $x, y \in{\mathbb R}^J$ and a closed set $A \subset{\mathbb R}^J$,
$d(x,y)$ denotes the Euclidean distance between $x$ and~$y$, and
$d(x,A) = \inf_{y \in A} d(x,y)$ denotes the distance between $x$ and
the set $A$. For each $r\geq0$, $N_r(A)=\{x\in{\mathbb R}^J \dvtx d(x,A)\leq r\}$.
The unit sphere in ${\mathbb R}^J$ is represented by $S_1(0)$.
Given a set $A \subset{\mathbb R}^J$, $A^\circ$ denotes its interior,
$\overline{A}$ its closure and
$\partial A$ its boundary.

The space of
continuous functions on $[0,\infty)$ that take values in ${\mathbb
R}^J$ is
denoted by ${\mathcaligr C}[0,\infty)$,
and,
given a set $G \subset{\mathbb R}^J$,
${\mathcaligr C}_G [0,\infty)$ denotes the subset of functions~$f$ in
${\mathcaligr C}[0,\infty)$ such that
$f(0) \in G$. The spaces ${\mathcaligr C}[0,\infty)$ and ${\mathcaligr
C}_G [0,\infty)$ are assumed to be equipped
with the topology of uniform convergence on compact sets.
Given $f \in{\mathcaligr C}[0,\infty)$ and $T \in[0,\infty)$,
$\mathrm{Var}_{[0,T]} f$
denotes the
${\mathbb R}_+ \cup\{\infty\}$-valued number
that equals the variation of $f$ on $[0,T]$.
Also, given a real-valued function $f$ on $[0,\infty)$, its
oscillation is defined by
\[
\mathit{Osc}(f; [s,t]) = \sup_{s\leq u_1 \leq u_2 \leq t} |f(u_2) - f(u_1)|;\qquad
 0\leq s \leq t < \infty.
\]
For each $A\in{\mathbb R}^J$, $\mathbb{I}_A(\cdot)$ denotes the
indicator function
of the
set $A$, which takes the value $1$ on $A$ and $0$ on the complement of $A$.

Given two random variables $U^{(i)}$ defined on a probability\vspace*{-1pt} space
 $(\Omega^{(i)}, {\mathcaligr F}^{(i)},\break {\mathbb P}^{(i)})$
and taking
values in a common Polish space\vspace*{1pt} $S$, $i = 1, 2$,
the notation $U^{(1)} \stackrel{(d)}{=} U^{(2)}$ will be used to imply
that the random variables are equal in
distribution. Given a sequence of $S$-valued random variables $\{
U^{(n)}, n \in{\mathbb N}\}$ and $U$, with $U^{(n)}$ defined on
$(\Omega^{(n)}, {\mathcaligr F}^{(n)}, {\mathbb P}^{(n)})$ and $U$
defined on
$(\Omega, {\mathcaligr F}, {\mathbb P})$,
$U^{(n)} \Rightarrow U$ is used to denote weak convergence of the
sequence $U^{(n)}$ to $U$.
Also, if the sequence of random variables are all defined on the same
probability space $(\Omega, {\mathcaligr F}, {\mathbb P})$,
the notation $U^{(n)} \stackrel{({\mathbb P})}{\rightarrow} 0$ is
used to denote
convergence in probability.

\section{The class of reflected diffusions}
\label{sec-mainres}

The class of stochastic differential equations with reflection under
study are introduced in Section \ref{subs-sder}, and the basic
assumptions are stated in Section \ref{subs-ass}.
Some useful ramifications of the assumptions and a motivating example
are then presented in Section \ref{subs-motiv}.

\subsection{Stochastic differential equations with reflection}
\label{subs-sder}

The so-called extended Skorokhod problem (ESP), introduced in \cite{ram1},
is a convenient tool for the pathwise construction of reflected diffusions.
The data associated with an ESP is the closure~$G$ of an open,
connected domain in ${\mathbb R}^J$ and
a set-valued mapping $d(\cdot)$ defined on $G$ such that $d(x) = \{0\}
$ for $x \in G^\circ$,
$d(x)$ is a nonempty, closed and convex cone in ${\mathbb R}^J$ with
vertex at
the origin for every
$x \in\partial G$ and the graph of $d(\cdot)$ is closed.
Roughly speaking, given a continuous path $\psi$, the ESP associated
with $(G,d(\cdot))$
produces a constrained version
$\phi$ of $\psi$ that is restricted to live within the domain $G$ by
adding to it a ``constraining term''
$\eta$ whose increments over any interval lie in the closure of the
convex hull of the union of the allowable directions $d(x)$
at the points $x$
visited by $\phi$ during this interval.
We now state the rigorous definition of the ESP.
(In \cite{ram1}, the ESP was formulated more generally
for c\`{a}dl\`{a}g paths, but the formulation below will
suffice for our purposes since we consider only
continuous processes.)

\begin{defn} [(Extended Skorokhod problem)]
\label{def-esp}
Suppose $(G,d(\cdot))$ and $\psi\in{\mathcaligr C}_G [0,\infty)$
are given.
Then $(\phi, \eta) \in{\mathcaligr C}_G [0,\infty)\times
{\mathcaligr C}[0,\infty)$
are said to solve the ESP for $\psi$
if $\phi(0) = \psi(0)$,
and if for all
$t \in[0, \infty)$, the following properties hold:
\begin{enumerate}
\item
$\phi(t) = \psi(t) + \eta(t)$;
\item
$\phi(t) \in G$;
\item
for every $s \in[0, t)$
%
\begin{equation}
\label{hullprop}
\eta(t) - \eta(s) \in\overline{\operatorname{co}}  \biggl[ \bigcup_{u \in(s,t]} d(\phi
(u))   \biggr],
\end{equation}
where $\overline{\operatorname{co}}[A]$ represents the closure of the
convex hull
generated by the
set $A$.
\end{enumerate}
If $(\phi,\eta)$ is the unique solution to the ESP for $\psi$,
then we write $\phi= \bar{\Gamma}(\psi)$, and refer to
$\bar{\Gamma}$ as the extended Skorokhod map (ESM).
\end{defn}

If a unique solution to the ESP exists for all $\psi\in{\mathcaligr
C}_G [0,\infty)$, then
the associated ESM $\bar{\Gamma}$ is said to be well-defined on
${\mathcaligr C}_G [0,\infty)$.
In this case, it is easily verified (see Lemma \ref{Lip}) that
if $\phi=\bar{\Gamma}(\psi)$, then for any $s \in[0,\infty)$,
$\phi^s = \bar{\Gamma}(\psi^s)$, where for $t \in[0,\infty)$,
%
\begin{equation}
\label{shift}
\psi^s (t) \doteq\phi(s) + \psi(s + t) - \psi(s),\qquad  \phi^s
(t) \doteq\phi(s + t).
\end{equation}
Moreover, a well-defined
ESM is said to be Lipschitz continuous on
${\mathcaligr C}_G [0,\infty)$ if for every $T < \infty$, there exists
$\overline{K}_T < \infty$ such that,
given $\psi^{(i)} \in{\mathcaligr C}_G [0,\infty)$
with corresponding solution $(\phi^{(i)}, \eta^{(i)})$ to the ESP,
for $i = 1,
2$, we have
%
\begin{equation}
\label{ineq-lip}
\sup_{s \in[0,T]} \big|\phi^{(1)} (s) - \phi^{(2)}(s) \big|
\leq\overline{K}_T \sup_{s \in[0,T]} \big|\psi^{(1)} (s) - \psi^{(2)}
(s) \big|.
\end{equation}

The ESP is a generalization of the so-called Skorokhod Problem (SP)
introduced in \cite{sko2}.
Unlike the SP, the ESP does not require that the constraining term
$\eta$
have finite variation on bounded intervals
(compare Definitions 1.1 and 1.2 of \cite{ram1}).
The ESP can be used to define solutions to
stochastic differential equations with reflection (SDERs) associated
with a given
pair $(G,d(\cdot))$ and drift and dispersion coefficients $b\dvtx {\mathbb R}^J
\mapsto{\mathbb R}^J$ and
$\sigma\dvtx {\mathbb R}^J \mapsto{\mathbb R}^J \times{\mathbb R}^N$.

\begin{defn}
\label{def-SDER}
Given $(G,d(\cdot))$, $b(\cdot)$ and $\sigma(\cdot)$, the
triple
$(Z_t, B_t), (\Omega, {\mathcaligr F},\break {\mathbb P}), \{{\mathcaligr
F}_t\}$
is said to be a weak solution to the associated SDER if and only if:
\begin{enumerate}
\item
$\{{\mathcaligr F}_t\}$ is a filtration on the probability space
$(\Omega, {\mathcaligr F}, {\mathbb P})$ that satisfies the usual conditions;
\item
$\{B_t, {\mathcaligr F}_t\}$ is an $N$-dimensional Brownian motion.
\item
${\mathbb P}  ( \int_0^t   | b(Z(s))   | \, ds +
\int_0^t   | \sigma(Z(s))   |^2 \, ds < \infty  ) = 1$
 $\forall t \in[0,\infty)$.
\item
$\{Z_t, {\mathcaligr F}_t\}$ is a $J$-dimensional, adapted process
such that
${\mathbb P}$-a.s., $(Z, Y)$ solves the ESP for $X$, where $Y \doteq Z
- X$ and
%
\begin{equation}
\label{def-x}
  X(t) = Z(0) + \int_0^t b(Z(s)) \, ds + \int_0^t \sigma(Z(s))
\, d B(s)\qquad  \forall t \in[0,\infty).
\end{equation}
\item
The set $\{t \dvtx Z(t) \in\partial G\}$ has ${\mathbb P}$-a.s. zero Lebesgue
measure. In
other words, \mbox{${\mathbb P}$-a.s.},
%
\begin{equation}
\label{eq-zero}
\int_0^\infty\mathbb{I}_{{\partial G}} (Z(s)) \, ds = 0.
\end{equation}
\end{enumerate}
%
\end{defn}

This is similar to the usual definition for weak solutions for SDEs
(see, e.g.,
Definitions 3.1 and 3.2 of \cite{karshrbook}), except that property 4
is modified to
define reflection and property 5 captures the notion of ``instantaneous''
reflection (see, e.g., pages~87--88 of \cite{freibook}).
A strong solution can also be defined in an analogous
fashion.

\begin{defn}
\label{def-SDERstrong}
Given a probability space $(\Omega, {\mathcaligr F}, {\mathbb P})$ and
an $N$-dimensional Brownian motion $B$ on $(\Omega, {\mathcaligr F},
{\mathbb P})$,
$Z$ is said to be a strong solution to the SDER associated with
$(G,d(\cdot))$, $b(\cdot)$, $\sigma(\cdot)$ and
initial condition $\xi$ if
${\mathbb P}(Z(0) = \xi) =1$ and properties 3--5 of Definition \ref{def-SDER}
hold with $\{{\mathcaligr F}_t\}$ equal to the
completed and augmented filtration generated by the Brownian motion $B$.
\end{defn}

For a precise construction of the filtration $\{{\mathcaligr F}_t\}$ referred
to in
Definition \ref{def-SDERstrong}, see~(2.3) of \cite{karshrbook}.
In what follows, given the constraining process $Y$ in property 4 of
Definition \ref{def-SDER},
the quantity $L$ will denote the associated total variation measure:
in other words, for $0 \leq s \leq t < \infty$, we define
%
\begin{equation}
\label{eq-l}
L(s,t)\doteq{\mathrm{Var}}_{(s,t]}Y  \quad\mbox{and}\quad   L(t)
\doteq L(0,t].
\end{equation}
Observe that the process $L$ in the second definition in (\ref{eq-l}) is
$\{{\mathcaligr F}_t\}$-adapted and
takes values in the extended nonnegative reals, $\overline{{\mathbb R}}_+$.

\subsection{Main assumptions}
\label{subs-ass}

We now introduce certain basic assumptions on $(G,d(\cdot))$, $b(\cdot
)$ and
$\sigma(\cdot)$ that will be used in this work.
In Section \ref{subs-motiv},
we provide a concrete motivating example of a family of SDERs
that arise in applications which satisfies all the stated assumptions.
In Section \ref{sec-egburtob}, we provide another example of a class
of SDERs
that satisfy these assumptions. The latter class, which consists of
two-dimensional reflected diffusions in curved domains,
was first studied by Burdzy and Toby in \cite{burtoby}.

The first assumption concerns existence of solutions.
General conditions on $G$ and $d(\cdot)$ under
which this assumption is satisfied can be found in
Lemma 2.6, Theorem 3.3 and Theorem 4.3 of \cite{ram1}.

\begin{ass}
\label{as-weak}
There exists a weak solution
$(Z_t, B_t), (\Omega, {\mathcaligr F}, {\mathbb P}), \{{\mathcaligr
F}_t\}$
to the SDER associated with
$(G,d(\cdot))$, $b(\cdot)$ and $\sigma(\cdot)$ such that
$\{Z_t, {\mathcaligr F}_t; t \geq0\}$ is a Markov process under
${\mathbb P}$.
\end{ass}

Next, we impose a kind of ${{\mathbb L}}^p$-continuity condition on the ESM.

\begin{ass}
\label{as-weakcont}
There exist $p > 1, q \geq2$ and $K_T < \infty$, $T \in(0,\infty)$,
such that the weak solution $Z$ to the
SDER satisfies, for every $0 \leq s \leq t \leq T$,
%
\begin{equation}
\label{l2cont}
{\mathbb E}  [   | Y(t) - Y(s)   |^p \vert {\mathcaligr F}_s
  ] \leq
K_T {\mathbb E}  \Bigl[\sup_{u \in[s,t]}   |X(u) - X(s)  |^q
\big\vert
{\mathcaligr F}_s   \Bigr],
\end{equation}
where $X$ is the process defined by (\ref{def-x}) and $Y \doteq Z - X$.
\end{ass}

\begin{remark}
\label{rem-weakcont}
 Assumption \ref{as-weakcont} holds under rather mild conditions
on the ESP---for example,
when the following oscillation inequality
is satisfied for any solution $(\phi, \eta)$ to the ESP for a given
$\psi$: for every $0 \leq s \leq t < \infty$, there
exists $C_{s,t} < \infty$ such that
\[
\mathit{Osc} (\phi, [s,t]) \leq C_{s,t} \mathit{Osc} (\psi, [s,t]).
\]
In this case, since
$(Z, Y)$ solve the ESP for $X$, we have
for $0 \leq s \leq t \leq T$,
\[
|Y(t) - Y(s)| \leq \mathit{Osc} (Y, [s,t]) \leq C_{s,t} \mathit{Osc} (X, [s,t])
\leq2 C_{T} \sup_{u \in[s,t]} |X(u) - X(s)|,
\]
where $C_T = \max_{0 \leq s \leq t \leq T} C_{s,t} < \infty$,
and so Assumption \ref{as-weakcont} holds with $p = q = 2$ and $K_T =
4 C_{T}^2$.
The oscillation inequality can be shown to hold in many situations of interest
(see, e.g., Lemma 2.1 of \cite{wilsur95}).
If the ESM associated with $(G,d(\cdot))$ is well-defined and
Lipschitz continuous on ${\mathcaligr C}_G [0,\infty)$,
then the oscillation inequality is also automatically satisfied, and so
Assumption
\ref{as-weakcont} again holds with $p = q = 2$. Furthermore,
it is easy to see that if the ESM is well-defined and
H\"{o}lder continuous on ${\mathcaligr C}_G [0,\infty)$
with some exponent $\alpha\in(0, 1)$, then
Assumption \ref{as-weakcont} holds for any $p \geq2/\alpha$
and $q = \alpha p$. An example of such an ESM is provided in
Section \ref{sec-egburtob} (see also Section \ref{subs-pfegburtob}
and, in
particular, Remark \ref{rem-egburtob}).
\end{remark}

\begin{ass}
\label{as-locbd}
The coefficients $b$ and $\sigma$ are locally bounded, that is, they are
bounded on every compact subset of $G$.
\end{ass}

\subsection{A motivating example and ramifications of the assumptions}
\label{subs-motiv}

We now describe a family of multi-dimensional
ESPs $(G,d(\cdot))$ that arise in applications.
Fix $J \in{\mathbb N}$, $J \geq2$.
The $J$-dimensional ESPs in this family have domain $G = {\mathbb
R}_+^J$ and a
constraint
vector field $d(\cdot)$ that is parametrized by a ``weight'' vector
$\alpha= (\alpha_1, \ldots, \alpha_J)$ with $\alpha_i > 0$, $i =1,
\ldots,
J$, and $\sum_{i=1}^J \alpha_i =1$.
Associated with each weight vector $\alpha$
is the ESP $({\mathbb R}_+^J, d(\cdot))$, where
for $x \in\partial G = \partial{\mathbb R}_+^J$,
\[
d(x) \doteq  \biggl\{ \sum_{i:x_i = 0} \beta_i d_i \dvtx \beta_i \geq
0  \biggr\}
\]
with
\[
(d_{i})_j \doteq
\cases{\displaystyle
-\frac{\alpha_j}{1 - \alpha_i} , &\quad  for  $j \neq i$,\vspace*{2pt} \cr
1, &\quad for $j = i$,
}
\]
for $i, j = 1, \ldots, J$.
Reflected diffusions associated with this family were shown in \cite{ramrei1,ramrei2} to arise
as heavy traffic approximations of the so-called generalized processor
sharing (GPS) model in communication networks
(see also \cite{dupram98} and \cite{dupram99b}).
Indeed, the characterization of this class of reflected diffusions
serves as one of the motivations for this work.

Next, we introduce a family of SDERs that is a slight generalization
of the class of GPS ESPs.

\begin{defn}
\label{def-gpssder}
We will say $(G,d(\cdot)), b(\cdot)$ and $\sigma(\cdot)$ define a
{Class ${\mathcaligr A}$}
SDER if they satisfy the following conditions:
\begin{enumerate}
\item
The ESM associated with the ESP $(G,d(\cdot))$ is well-defined and,
for every
$T < \infty$, is Lipschitz continuous (with constant $K_T < \infty$) on
${\mathcaligr C}_G[0,T]$.
\item
$G$ is a closed convex cone with vertex at the origin, ${\mathcaligr
V} = \{0\}
$ and there exists
$\vec{{\mathbf v}}\in G$ such that
\[
\langle\vec{{\mathbf v}}, d \rangle= 0 \qquad \mbox{for all } d
\in d(x)
\cap S_1(0),
 x \in\partial G\setminus\{0\};
\]
\item
\label{as-lc}
There exists a constant $\lcconst< \infty$ such that for all $x, y
\in G$,
\[
| \sigma(x) - \sigma(y) | + |b(x) - b(y)| \leq\lcconst|x - y|
\]
and
\[
|\sigma(x)| \leq\lcconst ,\qquad |b(x)| \leq\lcconst(1 +
|x|).
\]
\item
The covariance function $a\dvtx G \rightarrow{\mathbb R}^J \times{\mathbb R}^J$
defined by
$a(\cdot) = \sigma^T (\cdot) \sigma(\cdot)$ is uniformly
elliptic, that is, there exists $\lambda> 0$ such that
%
\begin{equation}
\label{eq-ue}
u^T a(x) u \geq\lambda|u|^2  \qquad\mbox{for all }
u \in{\mathbb R}^J, x \in G.
\end{equation}
\end{enumerate}
\end{defn}

We expect that the conditions in property 3 can be relaxed to a
local Lipschitz and linear growth condition on both
$b$ and $\sigma$, and the main result can still be proved by using
localization along with the current arguments.
However, to keep the notation simple, we impose the slightly stronger
assumptions above.

\begin{remark}
\label{rem-gps}
 ESPs in the GPS family defined above were shown
to satisfy properties 1 and 2 (the latter with $\vec{{\mathbf v}}= e_1
+ \cdots
+ e_J$) of Definition \ref{def-gpssder}
in Theorem~3.6 and Lemma 3.4 of \cite{ram1}, respectively.
\end{remark}

In Theorem \ref{th-SDER}, we summarize some consequences of
Assumptions \ref{as-weak}--\ref{as-locbd}, and also
show that {Class ${\mathcaligr A}$} SDERs satisfy
these assumptions.
The proof essentially follows from Theorem 4.3 of \cite{ram1} and
Proposition 4.1 of \cite{kanram08c}.
The following set,
%
\begin{equation}
\label{def-calv}
{\mathcaligr V} \doteq\bigl\{x \in\partial G \dvtx \mbox{there exists } d \in S_1(0)
\mbox{
such that } \{d, -d\} \subset d(x) \bigr\},
\end{equation}
was shown in \cite{ram1} to play
an important role in the analysis.

\begin{theorem}
\label{th-SDER}
Suppose $(G,d(\cdot)), b(\cdot)$ and $\sigma(\cdot)$ satisfy
Assumptions \ref{as-weak} and~\ref{as-weakcont}, and let
$(Z_t, B_t), (\Omega, {\mathcaligr F}, {\mathbb P}), \{{\mathcaligr
F}_t\}$
be a weak
solution to the associated SDER.
Then~$Z$ is an ${\mathcaligr F}_t$-semimartingale on $[0,\stopzero)$, where
%
\begin{equation}
\label{def-tv}
\stopzero\doteq\inf\{t \geq0 \dvtx Z(t) \in{\mathcaligr V} \},
\end{equation}
and ${\mathbb P}$-a.s., $Z$ admits the decomposition
%
\begin{equation}
\label{eq-z}
Z(\cdot) = Z(0) + M(\cdot) + A(\cdot),
\end{equation}
where for $t \in[0,\stopzero)$,
%
\begin{equation}
\label{def-MA1}
M(t) \doteq\int_0^t \sigma(Z(s)) \cdot dB (s), \qquad   A(t)
\doteq\int_0^t b(Z(s)) \, ds + Y(t),
\end{equation}
and $Y$ has finite variation on $[0,t]$ and satisfies
%
\begin{equation}
\label{eq-y}
Y(t) = \int_0^t \gamma(s) \, dL(s),
\end{equation}
where $L$ is given by (\ref{eq-l}) and
$\gamma(s) \in d(Z(s))$, $dL$-a.e. $s \in[0,t]$.
Moreover, if $(G,d(\cdot))$, $b(\cdot)$ and $\sigma(\cdot)$ satisfy
properties
1 and 3 of Definition \ref{def-gpssder}, then they also satisfy
Assumption \ref{as-weak}, Assumption \ref{as-weakcont} (with $p = q = 2$)
and Assumption \ref{as-locbd}.
In this case,
$\{Z_t, {\mathcaligr F}_t\}$
is in fact the pathwise unique strong solution to the SDER,
is a strong Markov process and has
${\mathbb E}[|Z(t)|^2] < \infty$ for every $t \in(0,\infty)$ if
${\mathbb E}[|Z(0|^2] < \infty$.
\end{theorem}
\begin{pf}
Let $X$ be the process defined by (\ref{def-x}). Then $X$ is
clearly a semimartingale and property 4 of
Definition \ref{def-SDER} shows that ${\mathbb P}$-a.s., $(Z, Z-X)$ satisfy
the ESP for $X$. Moreover,
Theorem 2.9 of \cite{ram1} shows that $Y = Z- X$ has ${\mathbb
P}$-a.s. finite variation on any closed
sub-interval of $[0, T_{{\mathcaligr V}})$.
This shows that $Z$ is an $\mathcaligr F_t$-semimartingale on
$[0,\stopzero)$ with
the decomposition given in
(\ref{eq-z})--(\ref{eq-y}), and thus establishes the first assertion
of the theorem.

Next, suppose $(G,d(\cdot))$, $b(\cdot)$ and $\sigma(\cdot)$
satisfy properties 1 and 3 of Definition~\ref{def-gpssder}.
Then property 3 of Definition \ref{def-gpssder} implies
Assumption \ref{as-locbd} is satisfied. In addition,
by Remark \ref{rem-weakcont}, property 1 ensures that Assumption
\ref{as-weakcont} holds with $p=q \geq2$.
Moreover, Theorem 4.3 of \cite{ram1} and
Proposition 4.1 of \cite{kanram08c} show that, in fact,
the associated SDER admits a pathwise unique strong solution $Z$,
which is also a strong Markov process.
Thus, Assumption \ref{as-weak} is also satisfied.
Hence, we have shown that Assumptions
\ref{as-weak}--\ref{as-locbd} hold.
The last assertion of the theorem can be established using standard techniques,
by a modification of the proof in Theorem 4.3 of \cite{ram1},
in the same manner as this result is proved for strong solutions to SDEs,
and so we omit the details of the proof.
\end{pf}

We conclude this section by stating a consequence of property 2 of
Definition~\ref{def-gpssder} that will
be useful in the sequel.
Let $\Gamma_1$ denote the (extended) Skorokhod map associated with the
1-dimensional (extended) Skorokhod problem with $G = {\mathbb R}_+$ and $d(0)=
{\mathbb R}_+$, $d(x) = 0$ if $x > 0$. It is well known (see, e.g.,
\cite{sko2} or Lemma 3.6.14 of \cite{karshrbook}) that $\Gamma_1$ is
well-defined on ${\mathcaligr C}_{{\mathbb R}_+} [0,\infty)$, and in fact
has the explicit form
%
\begin{equation}
\label{def-gamma1}
\Gamma_1 (\psi) (t) = \psi(t) + \sup_{s \in[0,t]}   [ -\psi
(s)   ] \vee0.
\end{equation}

\begin{lemma}
\label{lem-vset}
Suppose that $(G,d(\cdot))$ satisfies property 2 of Definition \ref
{def-gpssder}.
If $(\phi,\eta)$ solves the associated ESP for
$\psi\in{\mathcaligr C}_G [0,\infty)$, then $\langle\phi, \vec
{{\mathbf v}}\rangle= \Gamma_1
(\langle\psi, \vec{{\mathbf v}}\rangle)$.
\end{lemma}

The proof of this lemma is exactly analogous to the proof of Corollary
3.5 of \cite{ram1},
and is thus omitted.

\subsection{Another motivating example}
\label{sec-egburtob}

We now describe a family of two-dimen\-sional reflecting Brownian motions
(henceforth abbreviated to RBMs)
in ``valley-shaped'' domains with vertex at the origin and
horizontal directions of reflection.
This family of reflected diffusions,
which was first studied in \cite{burtoby},
is parameterized by
two continuous real-valued functions $L$ and $R$ defined on
$[0,\infty)$, with $L(0) = R(0) = 0$ and $L(y) < R(y)$ for all
$y > 0$. The associated domain $G$ is then given by\looseness=-1
\[
G \doteq\{ (x,y) \in{\mathbb R}^2 \dvtx y \geq0, L(y) \leq x \leq R(y) \}.
\]
Let $\partial^1 G \doteq\{(x,y) \in\partial G \setminus(0,0)\dvtx
x = L(y) \}$ and, likewise, let $\partial^2 G
\doteq\{(x,y) \in\partial G \setminus(0,0)\dvtx x = R(y) \}$.
Then the reflection vector field is defined by
\begin{eqnarray*}
d(x,y)
= \cases{\displaystyle
(1,0), &\quad $(x,y) \in\partial^1 G$, \cr
(-1,0), &\quad $(x,y) \in\partial^2 G$, \cr
\{v\dvtx v_1 \geq0\}, &\quad $(x,y) = (0,0)$.
}
\end{eqnarray*}
Thus, there are two opposing, horizontal directions of reflection on
the two lateral boundaries, and then an additional
vertical reflection direction $(0,1)$ at $(0,0)$ to ensure that the
Brownian motion can be constrained within the domain.
To conform with the general structure of ESPs,
at $(0,0)$ we in fact define $d(\cdot)$ to be the convex cone
(which, in this case, equals a half-space) generated
by the three directions $(1,0)$, $(-1,0)$ and $(0,1)$.
Note that ${\mathcaligr V} = \{0\}$ for this ESM and,
when $L$ and $R$ are linear functions, this reflected diffusion
is a special case of the Class ${\mathcaligr A}$ SDER's introduced in
the last section.

It was shown in Theorem 1 of \cite{burtoby}
(see also Section 4.3 of \cite{burkanram08}) that the ESM~$\bar
{\Gamma}$
corresponding to $(G,d(\cdot))$ is well-defined and thus, when
$B$ is a standard two-dimensional Brownian motion,
$Z = \bar{\Gamma}(B)$ is a well-defined reflected Brownian motion
starting at $(0,0)$ and is also a Markov process (see Theorem 2 of
\cite{burtoby}).
In Proposition 4.13 of
\cite{burkanram08}, RBMs in this family were shown not to be semimartingales.
As an application of the results of this paper, we show that
when $L$ and $R$ are sufficiently regular, $Z$ nevertheless
admits a useful decomposition (see Corollary~\ref{cor-egburtob}).

\section{Statement of main results}
\label{subs-mainres}

Theorem \ref{th-SDER} shows that if ${\mathcaligr V} = \varnothing$
then $Z$
is a semimartingale.
In fact, it was shown in Theorem 1.3 of \cite{ram1} that when
${\mathcaligr
V} = \varnothing$,
the ESM coincides with the SM.
The main focus of this work
is to understand the behavior of reflected diffusions $Z$
associated with ESPs $(G,d(\cdot))$ for which ${\mathcaligr V} \neq
\varnothing
$, with the
GPS family being a representative example.
In \cite{ram1}, it was shown that for the GPS family of ESPs,
$Z$ is a semimartingale until the first time it hits the origin.
However, the first result of the present paper (Theorem \ref{th-main1})
shows that~$Z$ is not a semimartingale on $[0,\infty)$.

\begin{theorem}
\label{th-main1}
Suppose $(G,d(\cdot))$, $b(\cdot)$ and $\sigma(\cdot)$ describe a
{Class ${\mathcaligr A}$} SDER.
Then the unique pathwise solution $Z$ to the associated SDER is not a
semimartingale.
\end{theorem}

The proof of Theorem \ref{th-main1} is given in Section
\ref{subs-main1pf}. As mentioned in Section \ref{sec-intro}, for
the special case when $G$ is a convex wedge in ${\mathbb R}^2$ and the
directions of constraint on the two
faces are constant and point at each other, $b \equiv0$ and $\sigma$
is the identity matrix
(i.e., $Z$ is a reflected Brownian motion),
this result follows from Theorem 5 of \cite{wil85} (with the parameters
$\alpha= 1$ and the wedge angle $\pi$ less than $180^\circ$ therein).
The fact that, when $J=2$, the reflected Brownian motion $Z$ defined
here is the same as the reflected Brownian motion defined via the submartingale
formulation in \cite{wil85} follows from Theorem 1.4(2) of \cite{ram1}.
This two-dimensional result can also be viewed as a special case of
Proposition 4.13 of \cite{burkanram08}.
However, the proofs in \cite{burkanram08} and \cite{wil85} do not
seem to extend easily
to higher dimensions. In this paper, we take a different approach that
is applicable in arbitrary
dimensions and to more general diffusions, in particular providing a different
proof of the two-dimensional result mentioned above.

As is well known, when a process is a semimartingale,
${\mathcaligr C}^2$ functionals of the process can be characterized
using It\^
{o}'s formula.
Theorem \ref{th-main1} can thus be viewed as a somewhat negative
result since it suggests that
{Class ${\mathcaligr A}$} reflected diffusions and, in particular,
reflected diffusions associated
with the GPS family that arise in applications, may not possess
desirable properties.
However, we show in Corollary~\ref{cor-main2} that these diffusions
are indeed tractable by establishing that
they belong to the class of Dirichlet processes (in the sense of F\"{o}llmer).
This follows as a special case of a more general result, which is
stated below as Theorem \ref{th-main2}.

In order to state this result, we first
recall the definitions of zero $p$-variation processes and Dirichlet processes
(see, e.g., Theorem 2 of \cite{fol81b}).

\begin{defn}
For $p > 0$, a continuous process $A$ is of zero $p$-variation if and
only if for any $T>0$,
%
\begin{equation}\label{zero-pvar}
\sum_{t_i\in\pi^n}   | A(t_{i}) - A(t_{i-1})   |^p
\stackrel{({\mathbb P})}{\rightarrow} 0
\end{equation}
for any sequence $\{\pi^n\}$ of partitions of $[0,T]$ with
$\Delta(\pi^n) \doteq\max_{t_i\in\pi^n}(t_{i+1}-t_i)\rightarrow0$
as $n\rightarrow\infty$. If the process $A$ satisfies (\ref
{zero-pvar}) with $p = 2$, then $A$ is said to be
of zero energy.
\end{defn}

\begin{defn} \label{DP}
The stochastic process $Z$ is said to be a Dirichlet process if the
following decomposition holds:
%
\begin{equation}
\label{eq-decomp}
Z =M + A,
\end{equation}
where $M$ is an ${\mathcaligr F}_t$-adapted local martingale and $A$
is a continuous,
${\mathcaligr F}_t$-adapted, zero energy process with $A(0) = 0$.
\end{defn}

Note that this is weaker than the original definition of a Dirichlet
process given by F\"{o}llmer \cite{fol81b},
which requires that $M$ and $A$ in the decomposition (\ref{eq-decomp})
be square integrable and that $A$ satisfy ${\mathbb E}  [ \sum
_{t_i \in
\pi^n} |A_{t_i} - A_{t_{i-1}}|^2   ] \rightarrow0$ as
$\Delta(\pi^n) \rightarrow0$, rather than satisfy
(\ref{zero-pvar}) with $p = 2$.
However, our definition can be viewed as a localized version and
coincides with Definition 2.4 of \cite{coqjakmemslo06} (see also
Definition 12 of
\cite{rusval07}).

\begin{remark} \label{zeroEn} The decomposition of a Dirichlet
process $Z$, into a
local martingale and a zero energy process starting at $0$, is unique.
For any $p > 1$ and partition $\pi^n$ of $[0,T]$,
\[
\sum_{t_i\in\pi^n}|A(t_{i+1})-A(t_i)|^p \leq\max_{t_i \in\pi^n}
|A(t_{i+1}) - A (t_i)|^{p-1} {\mathrm{Var}}_{[0,T]}(A).
\]
Therefore, it follows that if $A$ is continuous and of finite variation,
then it is also of zero $p$-variation, for all $p > 1$.
In particular, this shows that the class of Dirichlet processes
generalizes the class of continuous semimartingales.
\end{remark}

\begin{theorem}
\label{th-main2}
Suppose $(G,d(\cdot))$, $b(\cdot)$ and $\sigma(\cdot)$ satisfy Assumptions
\ref{as-weak} and \ref{as-locbd}, let
$Z$ be an associated weak solution that satisfies Assumption \ref{as-weakcont}
for some $p > 1$,
and let $Y = Z - X$, where $X$ is defined by (\ref{def-x}).
Then $Y$ has zero $p$-variation.
\end{theorem}

%

As an immediate consequence of Theorem \ref{th-main2}, Definition \ref
{DP} and Theorem \ref{th-SDER}, we have the following result.

\begin{cor}
\label{cor-main2}
Suppose $(G,d(\cdot))$, $b(\cdot)$ and $\sigma(\cdot)$ satisfy
Assumptions \ref{as-weak} and~\ref{as-locbd},
and also Assumption \ref{as-weakcont} with $p = 2$. Then the
associated reflected diffusion is a Dirichlet process.
In particular, reflected diffusions associated with {Class~${\mathcaligr A}$} SDERs are
Dirichlet processes.
\end{cor}

The next consequence of Theorem \ref{th-main2} concerns the
class of reflected diffusions described in
Section \ref{sec-egburtob}.

\begin{cor}
\label{cor-egburtob}
Suppose that $L$ and $R$ are two continuous functions on
$[0,y]$ given by
%
\begin{equation}
\label{lrcond}
L(y) = -c_L y^{\alpha_L},\qquad   R(y) = c_R y^{\alpha_R},\qquad   y
\in[0,\infty),
\end{equation}
for some $\alpha_L, \alpha_R, c_L, c_R \in(0,\infty)$, and let
$\alpha= \min(\alpha_L, \alpha_R)$.
If $\alpha\geq1$, then the associated two-dimensional reflected diffusion $Z$ described in
Section \ref{sec-egburtob}
is a Dirichlet process, that is, admits the
decomposition $Z = B+A$, where $A$ is a process with zero quadratic
variation.
\end{cor}


It was shown in \cite{burkanram08} that, for every $\alpha>0$, $Z$ is not a
semimartingale. In contrast, Corollary \ref{cor-egburtob} establishes a positive result in this
direction, showing that even when the domain has a cusp-like shape
(i.e., corresponding to $\alpha > 1$), the reflected diffusion is a Dirichlet
process. This partially resolves the open question raised in \cite{burtoby}, as
mentioned in Section \ref{subs-pfegburtob}, when $L$ and $R$ are linear (i.e., when
$\alpha_L = \alpha_R = 1$) the domain is wedge-shaped and the reflected diffusion
$Z$ is associated with a Class $\mathcaligr{A}$ SDER. In this case, Corollary \ref{cor-egburtob}
follows from Corollary \ref{cor-main2}. The proof of Corollary \ref{cor-egburtob} in the general
case is given in Section \ref{subs-pfegburtob}. It is natural to expect that the reflected
diffusion would also be a Dirichlet process when $\alpha < 1$, since this
corresponds to nicer ``flatter'' domains. However, this does not
directly follow from the simple proof of Corollary \ref{cor-egburtob} given in Section
\ref{subs-pfegburtob} (see Remark~\ref{rem54}).

\section{Reflected diffusions associated with Class $\mathcaligr{A}$ SDERs}
\label{sec-smg}

Throughout this section, we will assume that $(G,d(\cdot))$, $b(\cdot
)$ and $\sigma(\cdot)$ describe a
{Class ${\mathcaligr A}$} SDER.
Let~$B$ be an $N$-dimensional Brownian motion on a given probability
space $(\Omega, {\mathcaligr F}, {\mathbb P})$, let $\{{\mathcaligr
F}_t\}$ be
the right-continuous augmentation of the filtration generated by $B$
(see Definition (2.3) given in
\cite{karshrbook}). Also,
let $Z$ be the pathwise unique strong solution to the associated SDER
(which exists by Theorem \ref{th-SDER}), let $X$ be defined by (\ref
{def-x}), let $Y \doteq Z - X$
and let $L$ be the total variation process of $Y$ as defined in (\ref{eq-l}).
We use ${\mathbb E}$ to denote expectation with respect to ${\mathbb
P}$ and, for $z
\in G$,
let~${\mathbb P}_z$ (resp., ${\mathbb E}_z$) denote the probability
(resp., expectation)
conditioned on $Z(0) = z$.

This section is devoted to the proof of Theorem \ref{th-main1}.
The key step is to show that the constraining process $Y$ in the
extended Skorokhod decomposition for $Z$ has ${\mathbb P}_0$-a.s. infinite variation.
More precisely, let $\vec{{\mathbf v}}$ be the vector that satisfies
property 3 of
Definition \ref{def-gpssder} and, for any given $\varepsilon\geq0$, consider
the hyperplane
%
\begin{equation}
\label{def-hyp}
H_\varepsilon\doteq  \{ x \in{\mathbb R}^d \dvtx \langle\vec
{{\mathbf v}}, x \rangle
= \varepsilon
\} \cap G,
\end{equation}
and let
%
\begin{equation}
\label{def-tauve}
\tau^\varepsilon\doteq\inf\{ t \geq0\dvtx Z(t) \in H_\varepsilon\}.
\end{equation}
We now state the key result in the proof of Theorem \ref{th-main1}.

\begin{theorem}
\label{th-mainstep}
There exists $T < \infty$ such that ${\mathbb P}_0(L(T) = \infty) > 0$.
\end{theorem}

A somewhat subtle point to note is that Theorem \ref{th-mainstep}
does not immediately establish the fact that $Z$ is not a
semimartingale because we do not know a priori that, if
$Z$ were a semimartingale, then its Doob decomposition
must be of the form $Z = M +A$ given in
(\ref{eq-z}) and (\ref{def-MA1}).
However, in Section \ref{subs-main1pf} (see Proposition~\ref{prop-smg})
we establish that this is indeed the case, thus
obtaining Theorem \ref{th-main1} from Theorem \ref{th-mainstep}.
First, in Section \ref{subs-zerodrift}, we establish Theorem \ref{th-mainstep}
for the case when $b \equiv0$.
The proof for the general case is obtained from this result via a
Girsanov transformation
in Section \ref{subs-gendrift}.

\subsection{The zero drift case}
\label{subs-zerodrift}

Throughout this section, we assume $b \equiv0$ and
establish the following result.

\begin{prop}
\label{prop-mainstep}
If $b \equiv0$, then
we have
%
\begin{equation}
\label{eq-expltauzero}
{\mathbb E}_0  \bigl[e^{-L(\tau^1)}  \bigr] = 0,
\end{equation}
and hence,
%
\begin{equation}
\label{eq-ltauinfty}
L(\tau^1) = \infty,\qquad {\mathbb P}_0\mbox{-a.s.}
\end{equation}
\end{prop}

When combined with Lemma \ref{lem-finite}, which shows that
${\mathbb P}_0(\tau^1 < \infty) = 1$ when $b \equiv0$, Proposition
\ref{prop-mainstep} yields Theorem \ref{th-mainstep}.
The proof of Proposition \ref{prop-mainstep} is given in
Section \ref{subsub-auxres}.
The proof relies on an upper bound for ${\mathbb E}_0[e^{-L(\tau
^1)}]$, which is
obtained in Section \ref{subsub-proofs}, and some weak convergence
results, which are established in Section \ref{subsub-weakconv}.

\subsubsection{An upper bound}
\label{subsub-proofs}

To begin with, we use
the strong Markov property of~$Z$ to
obtain an upper bound on ${\mathbb E}_0[e^{-L(\tau^1)}]$.
Recall the definition of $\tau^0$ given in (\ref{def-tauve})
with $\varepsilon= 0$, noting that $H_0 = \{0\}$ because
$G$ is a closed convex cone with vertex at $0$.
Moreover, for $\varepsilon> 0$, we recursively
define two sequences of random times $\{\tau^\varepsilon_n\}_{n \in
{\mathbb
N}}$ and
$\{\nstop^\varepsilon_n\}_{n \in{\mathbb N}}$ as follows: $\nstop
^\varepsilon_0
\doteq0$
and for
$n \in{\mathbb N}$,
%
\begin{eqnarray}
\label{def-tausigma}
\tau^{\varepsilon}_n
& \doteq& \inf  \{ t \geq\nstop
_{n-1}^\varepsilon \dvtx Z(t)
\in H_\varepsilon  \}, \nonumber\\[-8pt]\\[-8pt]
\nstop^{\varepsilon}_n
& \doteq& \inf  \{ t \geq\tau
_n^\varepsilon \dvtx
Z(t) \in H_0   \}.\nonumber
\end{eqnarray}
Since $Z$ is continuous and $H_\varepsilon$ and $H_0$ are closed, it is
clear that $\tau^0$,
$\tau_n^\varepsilon$ and $\nstop_n^\varepsilon$ are ${\mathcaligr
F}_t$-stopping times.
For conciseness, we will often denote $\tau^\varepsilon_1$ simply by
$\tau^\varepsilon$, since
this is consistent with the notation of $\tau^\varepsilon$ given in
(\ref{def-tauve}).

\begin{lemma}
\label{lem-ubound} For every $\varepsilon\in(0,1)$,
\[
{\mathbb E}_0   \bigl[e^{-L(\tau^1)}  \bigr] \leq
\frac{{\mathbb E}_0  [{\mathbb P}_{Z(\tau^\varepsilon)}   (
\tau^0
\geq\tau
^1  )   ]}{{\mathbb E}_0  [{\mathbb P}_{Z(\tau
^\varepsilon)}
  ( \tau^0
\geq\tau^1   )   ] + {\mathbb E}_0  [{\mathbb
E}_{Z(\tau^\varepsilon)}
[   (1 - e^{-L(\tau^0)}   )
\mathbb{I}_{\{\tau^0 < \tau^1\}}   ]  ]}.
\]
\end{lemma}
\begin{pf}
{F}rom the elementary inequality
\[
L(\tau^1) \geq\sum_{n=1}^{\infty} \bigl(L(\nstop^{\varepsilon}_n \wedge\tau^1) -
L(\tau^{\varepsilon}_n \wedge\tau^1)\bigr),
\]
it immediately follows that
%
\begin{equation}
\label{ineq-mainstep}
{\mathbb E}_0   \bigl[e^{-L(\tau^1)}   \bigr] \leq{\mathbb E}_0
  \bigl[e^{-\sum^{\infty}_{n=1}(L(\nstop
^{\varepsilon}_n
\wedge\tau^1)-L(\tau^{\varepsilon}_n \wedge\tau^1))}  \bigr].
\end{equation}
For $n \geq2$, $\alpha_n^\varepsilon\geq\alpha_1^\varepsilon$ and
$\tau_n^\varepsilon
\geq
\alpha_1^\varepsilon$. Hence, on
the set $\{\shortsig\geq\tau^1\}$, we have $\nstop_n^\varepsilon
\wedge
\tau^1 = \tau_n^\varepsilon\wedge\tau^1 = \tau^1$ for
every $n \geq2$.
Therefore, the right-hand side of (\ref{ineq-mainstep}) can be
decomposed as
\begin{eqnarray*}
{\mathbb E}_0
  \bigl[e^{-\sum^{\infty}_{n=1}(L(\nstop
^{\varepsilon}_n
\wedge\tau^1)-L(\tau^{\varepsilon}_n \wedge\tau^1))}  \bigr]
& = &
{\mathbb E}_0  \bigl[e^{-(L(\shortsig\wedge\tau^1) -
L(\tau^{\varepsilon} \wedge\tau^1))} \mathbb{I}_{\{\shortsig\geq
\tau
^1\}}  \bigr]\\
&&{} + {\mathbb E}_0
  \bigl[e^{-\sum^{\infty}_{n=1}(L(\nstop
^{\varepsilon}_n
\wedge\tau^1)-L(\tau^{\varepsilon}_n \wedge\tau^1))} \mathbb
{I}_{\{
\shortsig< \tau^1\}}  \bigr].
\end{eqnarray*}
Conditioning on ${\mathcaligr F}_{\shortsig}$, using the fact that
$\mathbb{I}_{\{
\shortsig< \tau^1\}}$,
$L(\shortsig\wedge\tau^1)$ and $L(\tau^\varepsilon\wedge\tau^1)$ are
${\mathcaligr
F}_{\shortsig}$-measurable, the strong Markov property of $Z$ and the fact
that $Z_1(\shortsig) = 0$,
last term above can be rewritten as
\begin{eqnarray*}
&&{\mathbb E}_0
  \bigl[e^{-\sum^{\infty}_{n=1}(L(\nstop
^{\varepsilon}_n
\wedge\tau^1)-L(\tau^{\varepsilon}_n \wedge\tau^1))} \mathbb
{I}_{\{
\shortsig< \tau^1\}}  \bigr] \\
 &&\qquad= {\mathbb E}_0
  \bigl[{\mathbb E}_0   \bigl[e^{-\sum^{\infty}_{n=1}
(L(\nstop^{\varepsilon}_n
\wedge\tau^1)-L(\tau^{\varepsilon}_n \wedge\tau^1))} \mathbb
{I}_{\{
\shortsig< \tau^1\}}  |
{\mathcaligr F}_{\shortsig}  \bigr]   \bigr]\\
  &&\qquad= {\mathbb E}_0
  \bigl[e^{-(L(\shortsig
\wedge\tau^1)-L(\tau^{\varepsilon} \wedge\tau^1))} \mathbb{I}_{\{
\shortsig< \tau^1\}}{\mathbb E}_0   \bigl[ e^{-\sum^{\infty}_{n=2}(L(\nstop^{\varepsilon}_n
\wedge\tau^1)-L(\tau^{\varepsilon}_n \wedge\tau^1))} |
{\mathcaligr
F}_{\shortsig}  \bigr]  \bigr] \\
  &&\qquad= {\mathbb E}_0
  \bigl[e^{-(L(\shortsig
\wedge\tau^1)-L(\tau^{\varepsilon} \wedge\tau^1))} \mathbb{I}_{\{
\shortsig< \tau^1\}}{\mathbb E}_{Z(\shortsig)}
 \bigl [ e^{-\sum^{\infty}_{n=1}(L(\nstop
^{\varepsilon}_n
\wedge\tau^1)-L(\tau^{\varepsilon}_n \wedge\tau^1))}   \bigr]
  \bigr] \\
  &&\qquad= {\mathbb E}_0
  \bigl[e^{-(L(\shortsig
\wedge\tau^1)-L(\tau^{\varepsilon} \wedge\tau^1))} \mathbb{I}_{\{
\shortsig< \tau^1\}}  \bigr] {\mathbb E}_{0}
  \bigl[ e^{-\sum^{\infty}_{n=1}(L(\nstop
^{\varepsilon}_n
\wedge\tau^1)-L(\tau^{\varepsilon}_n \wedge\tau^1))}   \bigr].
\end{eqnarray*}
Combining the last two assertions and rearranging terms, we obtain
\[
{\mathbb E}_0
 \bigl [e^{-\sum^{\infty}_{n=1}(L(\nstop
^{\varepsilon}_n
\wedge\tau^1)-L(\tau^{\varepsilon}_n \wedge\tau^1))}  \bigr]
= \frac{{\mathbb E}_0  [e^{-(L(\shortsig\wedge\tau^1) -
L(\tau^{\varepsilon} \wedge
\tau^1))}\mathbb{I}_{\{\shortsig\geq
\tau^1\}}  ]}{1-{\mathbb E}_0  [e^{-(L(\shortsig\wedge
\tau^1)-L(\tau^{\varepsilon} \wedge
\tau^1))}\mathbb{I}_{\{\shortsig< \tau^1\}}  ]}.
\]
Together with (\ref{ineq-mainstep}), this yields the inequality
%
\begin{equation}
\label{eq-ubound}
{\mathbb E}_0  \bigl [e^{-L(\tau^1)}  \bigr] \leq
\frac{{\mathbb E}_0  [e^{-(L(\shortsig\wedge\tau^1) -
L(\tau^{\varepsilon} \wedge
\tau^1))}\mathbb{I}_{\{\shortsig\geq
\tau^1\}}  ]}{1-{\mathbb E}_0  [e^{-(L(\shortsig\wedge
\tau^1)-L(\tau^{\varepsilon} \wedge
\tau^1))}\mathbb{I}_{\{\shortsig< \tau^1\}}  ]}.
\end{equation}

We now show that the upper bound stated in the lemma follows from
(\ref{eq-ubound}). Due to
the nonnegativity of $L(\shortsig\wedge\tau^1) -
L(\tau^{\varepsilon} \wedge\tau^1)$
and the strong Markov property of $Z$, we have
%
%
\begin{eqnarray}\label{eq-num1}
{\mathbb E}_0  \bigl[e^{-(L(\shortsig\wedge\tau^1) -
L(\tau^{\varepsilon} \wedge
\tau^1))}\mathbb{I}_{\{\shortsig\geq
\tau^1\}}  \bigr]
& \leq& {\mathbb E}_0  \bigl[ \mathbb{I}_{\{\shortsig\geq
\tau^1\}}  \bigr] \nonumber\\
& =& {\mathbb E}_0  \bigl[ {\mathbb E}_0  \bigl [ \mathbb{I}_{\{
\shortsig\geq
\tau^1\}} \vert {\mathcaligr F}_{\tau^\varepsilon}   \bigr]  \bigr] \\
& = & {\mathbb E}_0   \bigl[ {\mathbb P}_{Z(\tau^\varepsilon)}   (
\tau^0
\geq\tau^1
  )  \bigr],\nonumber
\end{eqnarray}
where recall that $\tau^0 = \inf\{ t \geq0 \dvtx Z(t) \in H_0\}$.
Similarly, once again conditioning on ${\mathcaligr F}_{\tau
^\varepsilon}$
and using
the strong Markov property of $Z$,
we obtain
\begin{eqnarray*}
&&{\mathbb E}_0   \bigl[e^{-(L(\shortsig\wedge\tau^1) -
L(\tau^{\varepsilon} \wedge\tau^1))}
\mathbb{I}_{\{\shortsig< \tau^1\}}  \bigr] \\
   &&\qquad= {\mathbb E}_0  \bigl[ {\mathbb E}_0  \bigl[
e^{-(L(\shortsig
\wedge\tau^1) -
L(\tau^{\varepsilon} \wedge\tau^1))}
\mathbb{I}_{\{\shortsig< \tau^1\}} | {\mathcaligr F}_{\tau
^\varepsilon}
  \bigr]
  \bigr] \\
   &&\qquad= {\mathbb E}_0  \bigl[{\mathbb E}_{Z(\tau
^{\varepsilon})}
\bigl[ e^{-L(\tau^0
\wedge\tau^1)} \mathbb{I}_{\{\tau^0 < \tau^1\}}  \bigr]  \bigr].
\end{eqnarray*}
Therefore,
%
\begin{eqnarray}
\label{eq-denom1}
&&1 - {\mathbb E}_0   \bigl[e^{-(L(\shortsig\wedge\tau^1) -
L(\tau^{\varepsilon} \wedge\tau^1))}
\mathbb{I}_{\{\shortsig< \tau^1\}}  \bigr]\nonumber \\
&&\qquad =
{\mathbb E}_0  \bigl[1-{\mathbb E}_{Z(\tau^{\varepsilon})}   \bigl[
e^{-L(\tau^0
\wedge\tau^1)} \mathbb{I}_{\{\tau^0 < \tau^1\}}  \bigr]  \bigr] \\
&&\qquad ={\mathbb E}_0  \bigl[ {\mathbb P}_{Z(\tau^\varepsilon)} (\tau
^0 \geq
\tau^1)  \bigr]
+{\mathbb E}_0   \bigl[ {\mathbb E}_{Z(\tau^{\varepsilon})}
  \bigl[  \bigl(1-e^{-L(\tau^0)}  \bigr) \mathbb{I}_{\{\tau^0< \tau
^1\}}
  \bigr]  \bigr].\nonumber
\end{eqnarray}
The lemma follows from (\ref{eq-ubound}), (\ref{eq-num1}) and (\ref
{eq-denom1}).
\end{pf}

Next, we establish an elementary lemma that holds when the drift is zero.
Recall the vector $\vec{{\mathbf v}}$ of property 2 of Definition \ref
{def-gpssder}.

\begin{lemma}
\label{lem-00} When $b \equiv0$, the process
$\langle Z, \vec{{\mathbf v}}\rangle$ is an ${\mathcaligr
F}_t$-martingale on
$[0,\tau^0]$ and
for every $\varepsilon> 0$, ${\mathbb P}_0$-a.s.,
%
\begin{equation}
\label{eq-eps}
{\mathbb P}_{Z(\tau^\varepsilon)} (\tau^0 \geq\tau^1) =
\varepsilon.
\end{equation}
\end{lemma}
\begin{pf}
{F}irst, note that $ H_0 = \{0\} = {\mathcaligr V}$ by property 2 of
Definition \ref{def-gpssder} and so~$T_{\mathcaligr V}$
defined in (\ref{def-tv}) coincides with $\tau^0$. From
Lemma \ref{lem-vset} and the continuity of the sample paths of $Y$,
it follows that
for $t \in[0,\tau^0]$,
$\langle Y(t), \vec{{\mathbf v}}\rangle= 0$ and so ${\mathbb P}$-a.s.,
%
\begin{equation}
\label{eqn-zlemma}
\langle Z(t), \vec{{\mathbf v}}\rangle= \langle Z(0), \vec{{\mathbf
v}}\rangle+
\tilde{M},\qquad
  t \in[0, \tau^0],
\end{equation}
where $\tilde{M} \doteq\langle\int_0^\cdot\sigma(Z(s)) \cdot dB(s),
\vec{{\mathbf v}}\rangle$
is an ${\mathcaligr F}_t$ martingale on $[0,\tau^0]$ since $\sigma$
is uniformly
bounded by property 3 of Definition \ref{def-gpssder}.
This establishes the first assertion of the lemma.

The quadratic variation $\langle\tilde{M} \rangle$ of $\tilde{M}$ is
given by
\[
\langle\tilde{M} \rangle(t) = \int_{0}^{t} \vec{{\mathbf v}}^T a(Z(s))
\vec{{\mathbf v}}\, ds,\qquad   t \in[0,\infty),
\]
where $a \doteq\sigma^T \sigma$.
By property 4 of Definition \ref{def-gpssder}, $a(\cdot)$ is
uniformly elliptic.
Therefore, ${\mathbb P}$-a.s.,
$\langle\tilde{M} \rangle$ is strictly increasing and
$\langle\tilde{M} \rangle_{\infty} \doteq\lim_{t \rightarrow
\infty}
\langle
\tilde{M} \rangle(t) =\infty$. For $t \in[0,\infty)$, let
\[
T(t) \doteq\inf\{ s \geq0 \dvtx \langle\tilde{M}\rangle(s) > t\},\qquad
 {\mathcaligr G}_t \doteq{\mathcaligr F}_{T(t)},\qquad
\tilde{B}(t) \doteq\tilde{M}(T(t)).
\]
Then $\{\tilde{B}_t, {\mathcaligr G}_t\}_{t\geq0}$ is a standard
one-dimensional Brownian motion
(see, e.g., Theorem 4.6 on page 174 of \cite{karshrbook}).
Define $\tilde{\tau}^{\varepsilon} \doteq\inf\{ t \geq0 \dvtx \tilde
{B}(t) = \varepsilon\}$.
By (\ref{eqn-zlemma}), we have ${\mathbb P}_0$-a.s.,
\[
{\mathbb P}_{Z(\tau^{\varepsilon})}   ( \tau^0 \geq\tau
^1  ) =
{\mathbb P}  \bigl(\tilde{\tau}^0 \geq
\tilde{\tau}^1 \vert \tilde{B}(0) = \varepsilon  \bigr) =
\varepsilon,
\]
where the latter follows from well-known properties of Brownian motion.
This proves (\ref{eq-eps}).
\end{pf}

\begin{remark}
\label{rem-limmain}
From Lemmas \ref{lem-ubound} and \ref{lem-00}, we conclude that
for every $\varepsilon> 0$,
\[
{\mathbb E}  \bigl[e^{-L(\tau^1)}   \bigr] \leq\frac{\varepsilon
}{\varepsilon+
{\mathbb E}_{0}
  [ {\mathbb E}_{Z(\tau^\varepsilon)}
  [   ( 1 - e^{-L(\tau^\varepsilon)}   ) \mathbb{I}_{\{
\tau^0 <
\tau^1 \}}   ]   ]}.
\]
Thus, in order to establish (\ref{eq-expltauzero}),
it suffices to show that for some sequence $\{\varepsilon_k\}_{k \in
{\mathbb N}}$
such that
$\varepsilon_k \rightarrow0$ as $k \rightarrow\infty$,
\[
\operatorname{\lim\inf}\limits_{k \rightarrow\infty} \frac{1}{\varepsilon_k}
{\mathbb E}_{0}   \bigl[ {\mathbb E}_{Z(\tau^{\varepsilon_k})}
  \bigl[   \bigl( 1 - e^{-L(\tau^{\varepsilon_k})}   \bigr) \mathbb
{I}_{\{ \tau^0 <
\tau^1 \}}   \bigr]   \bigr] = \infty.
\]
This is established in Section \ref{subsub-auxres} using scaling
arguments. Since $Z$ is a reflected diffusion (rather than
just a reflected Brownian motion), the scaling arguments are
more involved and rely on some weak convergence results that are
established in Section \ref{subsub-weakconv}.
The reader may prefer to skip forward to the proof of Proposition~\ref
{prop-mainstep} in
Section \ref{subsub-auxres} and refer back to the results
in Section \ref{subsub-weakconv} when required.
\end{remark}

\subsubsection{A weak convergence result}
\label{subsub-weakconv}

Recall that we have assumed that the drift $b \equiv0$.
Now, let $\{\varepsilon_k\}_{k \in{\mathbb N}}$ and $\{x_k\}_{k \in
{\mathbb
N}}$ be
sequences such that
$\varepsilon_k \rightarrow0$ as $k \rightarrow\infty$ and $x_k \in
H_{\varepsilon
_k}$ for $k \in
{\mathbb N}$.
For each $k\in{\mathbb N}$, let $\newzk$ be the pathwise unique
solution to
the associated SDER with initial
condition $x_k$, and let $\newxk, \newyk$ and $\newlk$
be the associated processes as defined in Definition \ref{def-SDER} and
(\ref{eq-l}). For $k \in{\mathbb N}$, consider the scaled process
\[
B^k(t) \doteq\frac{B(\varepsilon_k^2 t)}{\varepsilon_k},\qquad   t
\in[0,\infty),
\]
which is a standard Brownian motion due to Brownian scaling.
Similarly, define
%
\begin{equation}
\label{def-scaling}
A^k (t) \doteq\frac{A_{(k)} ({\varepsilon_k}^2 t)}{\varepsilon_k},\qquad
  A = X, Y, Z, L,
\end{equation}
and let ${\mathcaligr F}^{k}_t \doteq{\mathcaligr F}_{{\varepsilon
_k}^2 t}$
for $t \in
[0,\infty)$.
Clearly, the processes
$Z^k$, $B^k$, $Y^k$ and $L^k$ are $\{{\mathcaligr F}^k_t\}$-adapted and
$L^k(t)={\mathrm{Var}}_{[0,t]}Y^k$ for every $t\geq0$. For
$(r, R) \in(0,\infty)^2$ such that $r<R$, let
%
\begin{equation}\label{def-thetaeps}\theta_{r, R}^{k} \doteq
\inf  \{t \geq0 \dvtx \langle Z^k (t), \vec{{\mathbf v}}\rangle\notin(r,
R)   \},\qquad   k \in{\mathbb N}.
\end{equation}

This section contains two main results. Roughly speaking, the first
result (Lemma \ref{lem-0})
shows that for the question under consideration, we can in effect
replace the state-dependent diffusion coefficient
$\sigma(\cdot)$ by $\sigma(0)$. This property is then used in
Corollary \ref{cor-0}
to provide bounds on the total variation sequence
$L^k(\theta_{r,R}^k)$, as $\varepsilon_k \rightarrow0$.
First, we observe that there exists a simple
equivalence between $(X^k, Z^k, Y^k)$ and another triplet of processes
that will be easier to work with.\looseness=-1

\begin{remark}
\label{rem-simple}
 For notational conciseness, we define the scaled diffusion coefficient
\[
\sigma^k (x) \doteq\sigma(\varepsilon^k x),\qquad   x \in{\mathbb
R}^J, k
\in{\mathbb N}.
\]
By the definition of $\newzk$ and the scaling (\ref{def-scaling}), it then
follows that
\[
X^k(t) = \frac{x_k}{\varepsilon_k} +\frac{1}{\varepsilon_k} \int
_0^{\varepsilon_k^2 t}
\sigma\bigl(Z_{(k)}(s)\bigr) \, dB(s) =
\frac{x_k}{\varepsilon_k} + \int_0^t \sigma^k (Z^k (s)) \, dB^k (s),
\]
where the last equality holds by the time-change theorem for stochastic
integrals
(see Proposition 1.4 in Chapter V of \cite{RY}).
This implies $Z^k$ is a strong solution to the
SDER associated with $(G, d(\cdot))$, $b \equiv0$, $\sigma^k$ and
the Brownian motion $\{B^k(t), {\mathcaligr F}^k_t\}_{t \geq0}$
defined on
$(\Omega, {\mathcaligr F}, {\mathbb P})$, with initial condition
$x_k/\varepsilon_k$.
If $\sigma$ satisfies properties 3 and 4 of Definition \ref
{def-gpssder} then so does $\sigma^k$, and thus
$(G,d(\cdot))$, $b \equiv0$ and~$\sigma^k$ also describe a {Class
${\mathcaligr A}$}
SDER. Therefore, by Theorem
\ref{th-SDER} there exists a pathwise unique solution $\tilde{Z}^k$
to the associated SDER for the Brownian motion
$\{B_t, {\mathcaligr F}_t\}$ with initial condition $x_k/\varepsilon_k$.
Let $\tilde{X}^k$ and $\tilde{Y}^k$ be the processes\vspace*{1pt} associated with
$\tilde{Z}^k$, defined in the usual manner as follows:
%
\begin{equation}
\label{def-tx}
\tilde{X}^k (t) = \frac{x_k}{\varepsilon_k} + \int_0^t \sigma^k
(\tilde
{Z}^k(s)) \, dB (s),\qquad   t \in[0,\infty),
\end{equation}
and $\tilde{Y}^k = \tilde{Z}^k - \tilde{X}^k$. From
the fact that solutions to {Class ${\mathcaligr A}$} SDERs are unique
in law by Theorem
\ref{th-SDER},
it then follows that
%
\begin{equation}
\label{claim-equiv}
(X^k, Z^k, Y^k) \stackrel{(d)}{=} (\tilde{X}^k, \tilde{Z}^k, \tilde{Y}^k),
\end{equation}
where recall that $\stackrel{(d)}{=}$ indicates equality in distribution.
\end{remark}

\begin{lemma}
\label{lem-0}
Given $x_* \in{\mathbb R}_+^J$, let
$(\overline{Z}, \overline{Y})$ satisfy the ESP pathwise for
%
\begin{equation}
\label{def-ox}
\overline{X} \doteq x_* + \sigma(0) B,
\end{equation}
and let
%
\begin{equation}
\label{def-thetalr}
\overline\theta_{r,R} \doteq\inf  \{t \geq0 \dvtx \langle
\overline Z(t), \vec{{\mathbf v}}\rangle\notin(r, R)   \}.
\end{equation}
Suppose $b \equiv0$ and $x_k/\varepsilon_k\rightarrow x_*$ as
$k\rightarrow
\infty$.
Then the following properties hold:
\begin{enumerate}
\item As $k \rightarrow\infty$,
%
\begin{equation}
\label{ZL2}
{\mathbb E}  \Bigl[\sup_{t\in[0,T]}| \tilde{Z}^k(t)-\overline
Z(t)|^2  \Bigr]
\rightarrow0
\end{equation}
and
$(X^k, Z^k, Y^k) \Rightarrow(\overline{X}, \overline{Z}, \overline{Y})$;
\item For
all but countably many pairs $(r, R) \in(0,\infty)^2$ such that
$r<R$, as $k \rightarrow\infty$, we have
\[
\max_{i=1, \ldots,J} \sup_{s \in[0,\theta_{r, R}^k]} Y^k_i(s)
\Rightarrow\max_{i=1, \ldots,J} \sup_{s \in[0,\overline{\theta
}_{r,
R}]} \overline{Y}_i(s).
\]
\end{enumerate}
\end{lemma}
\begin{pf} Note that since $x_k/\varepsilon_k \in H_1$ for every $k
\in
{\mathbb N}
$ and $H_1$ is closed, we must have $x_* \in H_1$.
We first prove property 1.
Let $\tilde{X}^k, \tilde{Z}^k$ and $\tilde{Y}^k$ be as in Remark
\ref{rem-simple}.
Then, by (\ref{claim-equiv}), it clearly suffices to show that
$(\tilde{X}^k, \tilde{Z}^k, \tilde{Y}^k) \Rightarrow
(\overline{X}, \overline{Z}, \overline{Y})$. From
(\ref{def-tx}) and (\ref{def-ox}), it follows that for $t \in
[0,\infty)$,
\begin{eqnarray*}
 |\tilde X^k(t)-\overline X(t)|^2
 &\leq&    \bigg|\frac{x_k}{\varepsilon_k} -x_* + \int
_0^t  \bigl(\sigma
^k(\tilde Z^k(s))-\sigma(0)  \bigr)\,dB(s)  \bigg|^2 \\
& \leq&    \biggl(  \bigg|\frac{x_k}{\varepsilon
_k}-x_*  \bigg|+  \bigg|\int
_0^t  \bigl(\sigma^k(\overline Z(s))-\sigma(0)  \bigr)\,dB(s) \bigg |
   \\
& &\hspace*{17.5pt}{} +    \bigg|\int_0^t  \biggl(\sigma
^k(\tilde
Z^k(s))-\sigma^k(\overline Z(s))  \biggr)\,dB(s) \bigg |   \biggr)^2.
\end{eqnarray*}
Using the fact that $(a+b+c)^2\leq3(a^2+b^2+c^2)$ for all $a,b,c\in
{\mathbb R}
$ and taking the supremum over $t \in[0,T]$ and then expectations of
both sides,
we obtain
\begin{eqnarray*}
 &&{\mathbb E}  \Bigl[\sup_{t\in[0,T]}|\tilde
X^k(t)-\overline
X(t)|^2
\Bigr]\\
&&\qquad \leq 3   \bigg|\frac{x_k}{\varepsilon_k}-
x_*  \bigg|^2 +3
{\mathbb E}
 \biggl [ \sup_{t\in[0,T]}
  \bigg|\int_0^t  \bigl(\sigma^k(\overline Z(s))-\sigma(0)
\bigr)\,dB(s)  \bigg|^2   \biggr] \\
 &&\quad\qquad{} +  3
{\mathbb E}  \biggl[\sup_{t\in[0,T]}  \bigg|\int_0^t  \bigl(\sigma
^k(\tilde
Z^k(s))-\sigma^k(\overline Z(s))  \bigr)\,dB(s) \bigg |^2  \biggr].
\end{eqnarray*}
Since $\sigma$ is uniformly bounded, the stochastic integrals on the
right-hand side are martingales.
By applying the Burkholder--Davis--Gundy (BDG) inequality, the
Lipschitz condition on $\sigma$, the definition of $\sigma^k$ and
Fubini's theorem, we obtain
\begin{eqnarray*}
&& {\mathbb E}  \biggl[\sup_{t\in[0,T]} \bigg |\int_0^t
\bigl(\sigma
^k(\tilde
Z^k(s))-\sigma^k(\overline Z(s))  \bigr)\,dB(s) \bigg |^2  \biggr] \\
 &&\qquad\leq  C_2{\mathbb E} \biggl [\int_0^T  \big|\sigma
^k(\tilde
Z^k(s))-\sigma^k(\overline Z(s))  \big|^2\,ds  \biggr] \\
  &&\qquad\leq  C_2\tilde K^2 \varepsilon_k^2 {\mathbb
E}  \biggl[\int_0^T
|\tilde
Z^k(s)-\overline Z(s)  |^2\,ds  \biggr] \\
 &&\qquad \leq  C_2\tilde K^2 \varepsilon_k^2 \int
_0^T{\mathbb E}  \Bigl[\sup
_{u\in
[0,s]}  |\tilde Z^k(u)-\overline Z(u)  |^2  \Bigr]\,ds,
\end{eqnarray*}
where $C_2 < \infty$ is the universal constant in the BDG inequality.
Using similar arguments, we also see that
\begin{eqnarray*}
 {\mathbb E}  \biggl[\sup_{t\in[0,T]} \bigg |\int_0^t
\bigl(\sigma
^k(\overline Z(s))-\sigma(0)  \bigr)\,dB(s) \bigg |^2  \biggr]
&\leq &C_2 \tilde{K}^2 \varepsilon_k^2 \int_0^T {\mathbb E} \Bigl [
\sup_{u
\in
[0,s]}   |\overline Z(u)  |^2  \Bigr] \, ds \\
& \leq&  C_2\tilde K^2\varepsilon_k^2 T
{\mathbb E}  \Bigl[\sup
_{t\in
[0,T]}  |\overline Z(t)  |^2  \Bigr].
\end{eqnarray*}
Combining the last three displays, and setting $\tilde{C}_T \doteq
3C_2\tilde{K}^2 (1 \vee T) < \infty$, we have
%
\begin{eqnarray}
\label{est:xk}
   {\mathbb E}  \Bigl[\sup_{t\in[0,T]}|\tilde
X^k(t)-\overline
X(t)|^2  \Bigr]
& \leq&  \tilde{C}_T \varepsilon_k^2
\int_0^T
{\mathbb E}  \Bigl[ \sup_{u \in[0,s]}   | \tilde{Z}^k (u) -
\overline
{Z}(u)   |^2   \Bigr] \, ds\nonumber\\[-8pt]\\[-8pt]
& & {}  + R^k(T),\nonumber
\end{eqnarray}
where
%
\[
R^k(T) \doteq3  \bigg| \frac{x_k}{\varepsilon_k} - x_*   \bigg|^2 +
\tilde
{C}_T \varepsilon_k^2
{\mathbb E}  \Bigl[\sup_{t\in[0,T]}  |\overline Z(t)
|^2  \Bigr].
\]
By the assumed Lipschitz continuity of $\overline{\Gamma}$,
\begin{eqnarray*}
{\mathbb E}  \Bigl[\sup_{t\in[0,T]}  |\overline Z(t)
|^2  \Bigr]
&\leq&
K_T^2
{\mathbb E} \Bigl [\sup_{t \in[0,T]}   |x_* + \sigma(0)
B(t)
|^2  \Bigr] \\
& \leq&
2 K_T^2 |x_*|^2 + 2 K_T^2 |\sigma(0)|^2
{\mathbb E}  \Bigl[\sup_{t \in[0,T]} |B(t)|^2   \Bigr]  < \infty.
\end{eqnarray*}
Since $x_k/\varepsilon_k \rightarrow x_*$ and $\varepsilon_k
\rightarrow0$ as $k \rightarrow
\infty$, it
follows that
%
\begin{equation}
\label{limit-rk}
\lim_{k \rightarrow\infty} R^k(T) = 0.
\end{equation}

On the other hand, combining the inequality in (\ref{est:xk}) with
the Lipschitz continuity of the map $\bar{\Gamma}$, we obtain
\begin{eqnarray*}
&&{\mathbb E}  \Bigl[\sup_{t\in[0,T]}|\tilde Z^k(t)-\overline
Z(t)|^2  \Bigr]\\
&&\qquad\leq K_T^2 R^k(T) + K_T^2 \tilde{C}_T \varepsilon_k^2 \int_0^T
{\mathbb E}  \Bigl[ \sup_{u \in[0,s]}   | \tilde{Z}^k (u) -
\overline
{Z}(u)   |^2   \Bigr] \, ds.
\end{eqnarray*}
An application of Gronwall's lemma then shows that
\[
{\mathbb E}  \Bigl[\sup_{t\in[0,T]}|\tilde Z^k(t)-\overline
Z(t)|^2  \Bigr]
\leq K_T^2
R^k(T) e^{K_T^2 \tilde{C}_T \varepsilon_k^2},
\]
which converges to zero as $k \rightarrow\infty$ due to
(\ref{limit-rk}) and the fact that $\varepsilon_k \rightarrow0$ as
$k \rightarrow
\infty$.
This proves
(\ref{ZL2}).
In turn, substituting the last inequality back into (\ref{est:xk})
and, again using (\ref{limit-rk}) and the
fact that $\varepsilon_k \rightarrow0$, we also obtain
\[
{\mathbb E} \Bigl [\sup_{t\in[0,T]}|\tilde{X}^k(t)-\overline
X(t)|^2  \Bigr]
\rightarrow0\qquad  \mbox{as } k \rightarrow\infty,
\]
which implies $\tilde{X}^k \Rightarrow\overline{X}$. Since the
mapping from $\tilde{X}^k \mapsto(\tilde{X}^k, \tilde{Z}^k, \tilde
{Y}^k)$ is continuous, by the continuous mapping theorem it follows that
$(\tilde{X}^k, \tilde{Z}^k, \tilde{Y}^k) \Rightarrow(\overline{X},
\overline{Z}, \overline{Y})$ and
the first property of the lemma is established.

We now turn to the proof of the second property.
By the first property, we know that $(Z^k, Y^k) \Rightarrow(\overline
{Z}, \overline{Y})$ as $k\rightarrow\infty$.
This immediately implies that for all but countably main pairs $(r,
R)
\in(0,\infty)^2$ such that $r<R$, we have, as $k \rightarrow\infty$,
\[
\bigl(Z^k(\cdot\wedge\theta^k_{r, R}),
Y^k (\cdot\wedge\theta^k_{r, R}), \theta^k_{r, R}\bigr)
\Rightarrow\bigl(\overline{Z} (\cdot\wedge\overline\theta_{r,
R}), \overline{Y}(\cdot\wedge\overline\theta_{r, R}),
\overline\theta_{r, R}\bigr).
\]
(For an argument that justifies this implication, see, e.g.,
the proof of Theorem 4.1 on page 354 of \cite{ethkurbook}.)
Using the continuity of the map
$(f, g, t) \mapsto\max_{i=1,\ldots, J} \sup_{s \in[0,t]} g_i(s)$ from
${\mathcaligr C}[0,\infty)\times{\mathcaligr C}[0,\infty)\times
{\mathbb R}_+$ to ${\mathbb R}_+$, an
application of
the continuous mapping theorem yields
the second property.
\end{pf}


\begin{cor}
\label{cor-0}
Suppose $b \equiv0$ and $x_k/\varepsilon_k\rightarrow x_*$ as
$k\rightarrow
\infty$. Then for each pair $(r, R) \in(0,\infty)$ such
that $r<R$, the following properties hold:
\begin{enumerate}
\item${\mathbb P}  (\sup_{k\in{\mathbb N}}L^k (\theta_{r,
R}^k)<\infty
  )=1$.
\item
$\varepsilon_k L^k (\theta_{r, R}^k) \Rightarrow0$.
\item
${\mathbb P}  (\overline{L}(\overline\theta_{r, R}) <
\infty
  ) =1$ and if $r<\langle  x,\vec{{\mathbf v}} \rangle <R$,
${\mathbb P}  (\overline{L}(\overline\theta_{r, R}) >
0  )
> 0$.
\end{enumerate}
\end{cor}
\begin{pf} If $\langle   x_*, \vec{{\mathbf v}} \rangle < r$ or $\langle   x_*,
\vec{{\mathbf v}} \rangle > R$, then $\overline\theta_{r, R}=0$
and $\theta_{r, R}^k=0$ for all $k$ sufficiently large. In
this case, properties (1)--(3) hold trivially. Hence, for the rest of the
proof, we assume that $r\leq\langle   x_*, \vec{{\mathbf v}}
\rangle\leq
R$.

We start by proving property 1. Let $\tilde{X}^k$, $\tilde{Z}^k$ and
$\tilde{Y}^k$ be defined as
in Remark \ref{rem-simple}, and let $\tilde{L}^k$ be defined as in
(\ref{eq-l}), but with $Y$ replaced by
$\tilde{Y}^k$. By (\ref{claim-equiv}), it follows that $(L^k, \theta
^k_{r, R})$ and $(\tilde L^k, \tilde{\theta}^k_{r, R
})$ have the same distribution for each $k\in{\mathbb N}$, where
$\tilde
{\theta}^k_{r, R}$ is defined in the obvious
way:
\[
\tilde\theta_{r, R}^{k} \doteq
\inf  \{t \geq0 \dvtx \langle\tilde Z^k (t), \vec{{\mathbf
v}}\rangle\notin
(r, R)   \}.
\]
We now argue that
${\mathbb P}  (\tilde{\theta}^k_{r, R} < \infty  )
= 1$.
Indeed, for $k \in{\mathbb N}$ such that $x_k /\varepsilon_k \notin(r,R)$
this holds
trivially. On the other hand, if $\tilde{Z}^k (0) = x_k/\varepsilon_k
\in
(r,R)$ then this follows
because Lemma \ref{lem-vset} and the uniform ellipticity condition show
that, on $(0, \tilde{\theta}_{r,R})$, $\langle\tilde{Z}^k(t),
\vec{{\mathbf v}}
\rangle=
\langle\tilde{X}^k(t), \vec{{\mathbf v}}\rangle$ is a continuous
martingale whose
quadratic variation
is strictly bounded away from zero. Thus, $\langle\tilde{Z}^k,
\vec{{\mathbf v}}
\rangle$
is ${\mathbb P}$-a.s. unbounded, and hence $\tilde{\theta}^k_{r,R}$ is
${\mathbb P}$-a.s. finite.
Therefore, to prove property 1, it suffices to show that
\[
{\mathbb P}  \Bigl(\sup_{k\in{\mathbb N}}\tilde L^k (\tilde\theta
_{r, R}^k
\wedge T) <\infty
  \Bigr)=1,\qquad
  T>0.
\]

Fix
$T\in(0, \infty)$.
Since $r>0$, there exists $\delta>0$ such that $\langle   y,
\vec{{\mathbf v}} \rangle <r$ for all $y$ with $|y|\leq\delta$.
Let $\tilde\kappa_{\delta}^k \doteq\inf\{t\geq0 \dvtx |\tilde
Z^k(t)|\leq\delta\}$.
Then $\theta_{r, R}^k\leq\tilde\kappa_\delta^k$ for all
$k\in{\mathbb N}$.
Let
\[
\tilde C^k \doteq\sup_{t\in[0,T]}|\tilde Z^k(t)|\vee|\tilde
X^k(t)|.
\]
By property 1 of Lemma \ref{lem-0}, it follows that $(\tilde X^k,
\tilde Z^k)
\Rightarrow(\overline X, \overline Z)$
as $k\rightarrow\infty$. Using the continuity of the map
$(f, g) \mapsto\sup_{s \in[0,T]} |f(s)|\vee|g(s)|$ from
${\mathcaligr C}[0,\infty)\times{\mathcaligr C}[0,\infty)$ to
${\mathbb R}_+$, an application of the
continuous mapping theorem yields
$\tilde C^k \Rightarrow\overline C$ as $k\rightarrow\infty$, where
$\overline C \doteq\sup_{t\in[0,T]}|\overline Z(t)|\vee|\overline
X(t)|$. Also, due to the Lipschitz continuity of the ESM
$\bar{\Gamma}$ and (\ref{def-ox}), ${\mathbb P}$-a.s., we have
\[
\sup_{t\in[0,T]}|\overline Z(t)| \leq K_T \sup_{t\in
[0,T]}|\overline X(t)|
\leq K_T  \Bigl(|x_*| + |\sigma(0)|\sup_{s \in[0,T]} |B(s)|  \Bigr) <
\infty,
\]
and hence ${\mathbb P}$-a.s.,
$\overline C <\infty$.
It then follows that ${\mathbb P}(\sup_{k\in{\mathbb N}}\tilde
C^k<\infty)=1$.
Moreover, ${\mathcaligr V} = \{0\}$ and for each $\omega\in\Omega$,
$(\overline Z(\cdot, \omega),\overline Y(\cdot,\omega))$ solves the
ESP for $\overline X(\cdot,\omega)$.
Therefore, it follows from Lemma 2.8 of \cite{ram1} that there exist
$\rho> 0$, independent of~$k$, a finite set ${\mathbb I}=\{1,\ldots,I\}$
and a
collection of open sets
$\{\mathcaligr O_i, i \in{\mathbb I}\}$ of ${\mathbb R}^J$ and associated
vectors $\{v_i
\in S_1(0), i\in{\mathbb I}\}$
that satisfy the following two properties:
\begin{enumerate}
\item$[\{x\in G\dvtx |x|\leq\overline C\}\setminus N_{\delta
/2}(0)^\circ]\subset[\bigcup_{i\in{\mathbb I}}\vspace*{1pt}\mathcaligr O_i]$.
\item If $y\in\{x\in G\dvtx |x|\leq\overline C\}\cap N_\rho(\mathcaligr
O_i)$ for some $i\in{\mathbb I}$ then
$
\langle d, v_i\rangle\geq\rho$ for every $d\in d(y)$
with $|d|=1$.
\end{enumerate}
Moreover, as in the proof of Theorem 2.9 of \cite{ram1},
for each $\omega\in\Omega$,
we can define a sequence $\{(\overline T_m(\omega),
\overline i_m(\omega)), m=0,1,\ldots\}$ defined recursively as follows. Let
$\overline T_0(\omega) \doteq0$ and let $\overline i_0(\omega) \in
{\mathbb I}$ be
such that $\overline Z(0,\omega) = x_* \in\mathcaligr O_{\overline
i_0(\omega)}$.
Note that because $x_* \in(r,R)$ implies
$|x| > \delta> \delta/2$, such an $i_0$ exists by property (1) above.
Next, for each $m=0,1,\ldots ,$ whenever
$\overline T_m(\omega) < \overline\kappa_{\delta/2}(\omega) \doteq
\inf\{t\geq0\dvtx \overline Z(t,\omega)\in N_{\delta/2}(0)\}$, define
\[
\overline T_{m+1}(\omega) \doteq\inf\bigl\{t>\overline T_{m}(\omega)\dvtx \overline
Z(t,\omega) \notin N_{\rho/2}\bigl(\mathcaligr O_{\overline i_m(\omega
)}\bigr)^\circ\mbox{
or } \overline Z(t,\omega) \in N_{\delta/2}(0)\bigr\}.
\]
If $\overline T_{m+1}(\omega)<T\wedge\overline\kappa_{\delta
/2}(\omega)$,
choose $\overline i_{m+1}(\omega) \in{\mathbb I}$ such that
$\overline
Z(\overline
T_{m+1}(\omega),\omega)\in\mathcaligr O_{\overline i_{m+1}(\omega)}$.
Note that
such
an $i_{m+1}(\omega)$ exists by property (1) above.\vspace*{-1pt}
Let $\overline{N}(\omega) < \infty$ be the smallest integer such
that $\overline T_{\overline{N}(\omega)}(\omega) \geq T\wedge
\overline\kappa_{\delta/2}(\omega)$ and redefine $\overline
T_{\overline{N}(\omega)}(\omega) = T\wedge\overline\kappa
_{\delta/2}(\omega)$. (Note that $\overline{N}(\omega)$ and $\{
(\overline T_m(\omega), \overline i_m(\omega)), m=0,1,\ldots\}$ are constructed
in the same way as $M$ and $\{T_m, m\in{\mathbb N}\}$ in Theorem 2.9
of \cite
{ram1}, except that we replace $\rho$
and $\delta$ by $\rho/2$ and $\delta/2$, respectively.)

Since, as shown in Lemma \ref{lem-0}, $(X^k, Z^k, Y^k) \Rightarrow
(\overline{X}, \overline{Z}, \overline{Y})$
as $k\rightarrow\infty$ and $(\overline{X}, \overline{Z}, \overline
{Y})$ has
continuous paths,
by invoking the Skorokhod representation theorem,
we may assume without loss of generality that there exists $\tilde
{\Omega}$
with
${\mathbb P}(\tilde{\Omega}) = 1$ such that for every $\omega\in
\tilde
{\Omega}$, $(X^k(\omega), Z^k(\omega),
Y^k(\omega)) \rightarrow(\overline{X}(\omega),\break \overline{Z}(\omega),
\overline{Y}(\omega))$ uniformly on $[0,T]$
as $k\rightarrow\infty$. Let $\bar k < \infty$ be such that for all
$k>\bar
k$, $\sup_{t \in[0,T]}|Z^k(t,\omega)-\overline Z(t,\omega)|<(\rho
\wedge
\delta)/4$. Then $Z^k(\cdot,\omega)$ will stay in $N_{\rho
}(\mathcaligr O_{\overline i_m(\omega)})$ during
the interval $[\overline T_m(\omega),\overline T_{m+1}(\omega))$.
Exactly as in the proof of Lemma 2.9 of \cite{ram1}
(note that the argument there only requires that $\phi(t)\in N_{\rho
}(\mathcaligr
O_{k_{m-1}})$ for $t\in[T_{m-1},T_m)$), we can then argue that
$\tilde L^k(T\wedge\tau_\delta^k(\omega),\omega)\leq(4\tilde
C^k(\omega)
\overline N(\omega))/\rho$ for $\omega\in\tilde{\Omega}$.
Together with
the fact that ${\mathbb P}(\sup_{k\in{\mathbb N}}\tilde C^k<\infty
)=1$ and
$\overline N(\omega)<\infty$ for each $\omega\in\Omega$, this
shows that
${\mathbb P}(\tilde L^k (\tilde\tau_{\delta}^k\wedge T) < \infty) = 1$.
Since $\tilde L^k (\tilde\theta_{r, R}^k\wedge T)\leq\tilde L^k
(\tilde\tau_{\delta}^k\wedge T)$, we then have ${\mathbb P}(\sup
_{k\in{\mathbb N}
}\tilde L^k
(\tilde\theta_{r, R}^k\wedge T)<\infty)=1$.
This completes the proof of property 1.

Property 2 follows directly from property 1 and the fact that
$\varepsilon
_k\rightarrow0$ as $k\rightarrow\infty$.
In addition, by Theorem \ref{th-SDER} it follows that $\overline Z$ is a
semimartingale on $[0,T_{\mathcaligr V})$, with~$\overline Y$ being
the bounded
variation term in the decomposition.
The first assertion of property 3 is thus a direct
consequence of the fact that $\overline\theta_{r, R
}<T_{\mathcaligr V}$.
For the second assertion of property 3, notice that with positive
probability, the Brownian motion
$\overline{X} = x_*+\sigma(0)B$ will exit $G$ before it hits one of
the two levels $H_{r}$
or $H_{R}$.
Since $\overline Z$ lies in $G$ and $\overline Z = \overline X +
\overline
Y$, this implies that, with positive probability, $\overline Y$ is not
identically zero in the interval $[0, \overline{\theta}_{r,R})$.
This, in turn, implies that $\overline L
(\overline{\theta}_{r,R})$ is strictly positive with positive probability.
Thus, the second assertion of property 3 is also established, and the
proof of
the corollary is complete.
\end{pf}

\subsubsection{A scaling argument}
\label{subsub-auxres}

Since the equality ${\mathbb E}_0[e^{-L(\tau^1)}] = 0$ implies that
${\mathbb P}$-a.s.,
$L(\tau^1) = \infty$,
in order to prove Proposition \ref{prop-mainstep} it suffices to establish
the former equality.
In turn, by Remark \ref{rem-limmain}, this equality holds if
there exists a sequence $\{\varepsilon_k\}_{k \in{\mathbb N}}$
such that $\varepsilon_k \rightarrow0$ as $k \rightarrow\infty$, and
%
\begin{equation}
\label{toshow-denom}
\operatorname{\lim\inf}\limits_{k \rightarrow\infty} \frac{1}{\varepsilon_k} {\mathbb
E}_0   \bigl[
{\mathbb E}_{Z(\tau
^{\varepsilon_k})}
  \bigl[  \bigl(1-e^{-L(\tau^0)}  \bigr) \mathbb{I}_{\{\tau^0< \tau
^1\}}
  \bigr]  \bigr] =\infty.
\end{equation}
We will show that (\ref{toshow-denom}) holds by using the strong
Markov property and scaling arguments.
First, we need to introduce some additional notation.
Fix $\varepsilon> 0$. Let $\Lambda_{\varepsilon}$ denote the
following union of hyperplanes:
%
\begin{equation}
\label{def-lam}
\Lambda_\varepsilon\doteq\bigcup_{n \in{\mathbb Z}} H_{2^n
\varepsilon}.
\end{equation}
For $x \in\Lambda_\varepsilon$, let $N_{\varepsilon}(x)$ denote the
pair of
hyperplanes in $\Lambda_\varepsilon$ that are
adjacent to the hyperplane on which $x$ lies. In other words, let
%
\begin{equation}
\label{def-neigh}
N_{\varepsilon} (x) \doteq H_{2^{n-1} \varepsilon} \cup H_{2^{n+1}
\varepsilon},\qquad
x \in H_{2^n \varepsilon},  n \in{\mathbb Z}.
\end{equation}
For future reference, note that for $y \in{\mathbb R}_+^J$ and
$x \in H_{2^n \varepsilon}$, $n \in{\mathbb Z}$,
%
\begin{equation}
\label{nscaling}
\frac{y}{\varepsilon} \in N_1   \biggl( \frac{x}{\varepsilon}
  \biggr)\quad
\Rightarrow\quad
  y \in N_{\varepsilon}
(x).
\end{equation}
Let $\{\newstop_n^\varepsilon\}_{n \in{\mathbb N}}$ be the
sequence of random times defined recursively
by $\newstop_{0}^\varepsilon\doteq0$ and for $n \in{\mathbb N}$,
%
\begin{equation}
\label{def-eta}
\newstop_{n}^\varepsilon\doteq\inf\{ t \geq\newstop
_{n-1}^\varepsilon\dvtx Z(t) \in
N_{\varepsilon} (Z(\newstop_{n-1}^\varepsilon)) \}.
\end{equation}
It is easy to see that $\{\newstop^\varepsilon_n\}_{n \in{\mathbb N}}$
defines a
sequence of stopping times
(for completeness, a proof is provided in Lemma \ref{lem-stop}).

Observe that $L$ is nondecreasing and for $x \in H_\varepsilon$,
${\mathbb P}
_x$-a.s., $\newstop^\varepsilon_n \leq\tau_0$
for every $n \in{\mathbb N}$. Now $Z(\tau^\varepsilon) \neq0$
because $\varepsilon>0$.
Hence, for every $n \in{\mathbb N}$,
%
\begin{equation}
\label{eq-denom2}
{\mathbb E}_{Z(\tau^{\varepsilon})}
  \bigl[  \bigl(1-e^{-L(\tau^0)}  \bigr) \mathbb{I}_{\{\tau^0< \tau
^1\}}
  \bigr] \geq{\mathbb E}_{Z(\tau^{\varepsilon})}
  \bigl[  \bigl(1-e^{-L(\newstop^\varepsilon_n)}  \bigr) \mathbb{I}_{\{
\tau^0< \tau
^1\}}   \bigr].
\end{equation}
Using the elementary identity
\[
1 -e^{-L(\newstop_n^\varepsilon)} = 1 - e^{-L(\newstop
_{n-1}^\varepsilon)} +
e^{-L(\newstop_{n-1}^\varepsilon)}   \bigl( 1 -
e^{-(L(\newstop_n^\varepsilon) - L(\newstop_{n-1}^\varepsilon))}
  \bigr),
\]
conditioning on ${\mathcaligr F}_{\newstop^\varepsilon_{n-1}}$, and
invoking the
strong Markov property of $Z$,
the right-hand side of (\ref{eq-denom2}) can be expanded as
\begin{eqnarray*}
&&{\mathbb E}_{Z(\tau^{\varepsilon})}
  \bigl[  \bigl(1-e^{-L(\newstop^\varepsilon_n)}  \bigr) \mathbb{I}_{\{
\tau^0< \tau
^1\}}   \bigr] \\
 &&\qquad= {\mathbb E}_{Z(\tau^{\varepsilon})} \bigl [
\bigl(1-e^{-L(\newstop
^\varepsilon_{n-1})}  \bigr) \mathbb{I}_{\{\tau^0< \tau^1\}}
\bigr] \\
 &&\quad\qquad{} + {\mathbb E}_{Z(\tau^{\varepsilon})} \bigl [ {\mathbb
E}_{Z(\tau
^{\varepsilon})}  \bigl[ e^{-L(\newstop^\varepsilon_{n-1})}  \bigl ( 1 -
e^{-(L(\newstop_n^\varepsilon) - L(\newstop_{n-1}^\varepsilon))}
  \bigr)\mathbb{I}_{\{
\tau^0< \tau^1\}}
 | {\mathcaligr F}_{\newstop^\varepsilon_{n-1}}   \bigr]   \bigr]
\\
&&\qquad = {\mathbb E}_{Z(\tau^{\varepsilon})}  \bigl[
\bigl(1-e^{-L(\newstop
^\varepsilon_{n-1})}  \bigr) \mathbb{I}_{\{\tau^0< \tau^1\}}
\bigr] \\
  &&\quad\qquad{}+ {\mathbb E}_{Z(\tau^{\varepsilon})}   \bigl[
e^{-L(\newstop
^\varepsilon_{n-1})}
{\mathbb E}_{Z(\newstop^\varepsilon_{n-1})}   \bigl[  \bigl ( 1 -
e^{-L(\newstop^\varepsilon
_1)}   \bigr) \mathbb{I}_{\{\tau^0 < \tau^1\}}   \bigr]   \bigr].
\end{eqnarray*}
Observing that the first term on the right-hand side is identical to
the term on the
left-hand side, except for a shift down in the index $n$,
we can iterate this procedure and use the relation $L(\newstop
_0^\varepsilon)
= L(0) = 0$
to conclude that for any $n \in{\mathbb N}$,
%
\begin{eqnarray}
\label{eq-denom3}
&&{\mathbb E}_{Z(\tau^{\varepsilon})}
  \bigl[ \bigl (1-e^{-L(\newstop^\varepsilon_n)}  \bigr) \mathbb{I}_{\{
\tau^0< \tau
^1\}}   \bigr] \nonumber\\[-8pt]\\[-8pt]
 && \qquad= \sum_{m=1}^{n} {\mathbb E}_{Z(\tau
^{\varepsilon})}
  \bigl[
e^{-L(\newstop^\varepsilon_{m-1})}
{\mathbb E}_{Z(\newstop^\varepsilon_{m-1})}   \bigl[   \bigl( 1 -
e^{-L(\newstop^\varepsilon
_1)}   \bigr) \mathbb{I}_{\{\tau^0 < \tau^1\}}   \bigr]   \bigr].\nonumber
\end{eqnarray}

Let $\{\varepsilon_k\}_{k \in{\mathbb N}}$ and $\{x_k\}_{k \in
{\mathbb N}}$
be sequences such
that $x_k \in H_{\varepsilon_k}$ for $k \in{\mathbb N}$, and
$\varepsilon_k \rightarrow0$ as $k \rightarrow\infty$. Since $H_1$
is compact and
$x_k/\varepsilon_k \in H_1$ for every $k \in{\mathbb N}$, we can
assume without
loss of generality (by choosing an appropriate subsequence, if
necessary) that there exists $x_* \in H_1$ such that $x_k/\varepsilon
_k \rightarrow x_*$,
as $k \rightarrow\infty$.
We now show that, when $\varepsilon$ is replaced by $\varepsilon_k$,
each term in the sum on the right-hand side of (\ref{eq-denom3}),
is $O(\varepsilon_k)$ (as $k \rightarrow\infty$),
with a constant that is independent of $m$.
This proof relies on the estimates obtained in the next two lemmas.
In both lemmas, $\newzk, \newyk, \newlk$, $Z^k$, $Y^k$ and $L^k$
denote the processes defined at the beginning
of Section \ref{subsub-weakconv},
and for $\varepsilon> 0$,
let $\newstop^{\varepsilon}_{(k),0} \doteq0$ and for $n \in{\mathbb N}$,
%
\begin{equation}
\label{def-newstopk1} \newstop^\varepsilon_{(k),n} \doteq\inf
\bigl\{ t
\geq\newstop^{\varepsilon}_{(k),n-1} \dvtx
\newzk(t) \in N_{\varepsilon}   \bigl(\newzk  \bigl(\newstop
^{\varepsilon}_{(k),n-1}
  \bigr)  \bigr)
  \bigr\},
\end{equation}
and, likewise, let $\newtime_{k,0} \doteq0$ and for $n \in{\mathbb
N}$, define
%
\begin{equation}
\label{def-newstopk2}
\newtime_{k,n} \doteq\inf  \bigl\{t \geq\newtime_{k,n-1} \dvtx
Z^k(t) \in N_{1}   \bigl( Z^k  (\newtime_{k,n-1}   )  \bigr)
  \bigr\}.
\end{equation}
Note that these sequences of stopping times are defined in a manner
analogous to the
sequence $\{\newstop_{n}^{\varepsilon}\}_{n \in{\mathbb N}}$
defined in (\ref{def-eta}), except that $Z$ is replaced by $\newzk$ and
$Z^k$, respectively.
Moreover, these definitions, together with the scaling relations
(\ref{def-scaling}) and (\ref{nscaling}),
yield the following equivalence relation
%
\begin{equation}
\label{newstop-equiv}
\varepsilon_k^2 \newtime_{k,n} = \beta_{(k),n}^{\varepsilon_k},\qquad
  k, n \in
{\mathbb N},
\end{equation}
\begin{lemma}
\label{lem-2}
Suppose $b \equiv0$. Then there exists $C > 0$ such that
%
\begin{equation}
\label{eq-lbound2}
\operatorname{\lim\inf}\limits_{k \rightarrow\infty} \frac{1}{\varepsilon_k} \inf_{x
\in
H_{\varepsilon_k}}
{\mathbb E}_{x}
 \bigl [   \bigl(1 - e^{-L (\newstop^{\varepsilon_k}_1)}   \bigr)
\mathbb{I}_{\{\tau
^0 <
\tau^1\}}   \bigr] \geq C.
\end{equation}
\end{lemma}

\begin{pf}
Since the law of $(\newzk,\newyk,L_{(k)})$ under ${\mathbb P}$ is the
same as
the law of $(Z,Y,L)$ under ${\mathbb P}_{x_k}$, we have
%
\begin{eqnarray}
\label{equivP}
 &&\operatorname{\lim\inf}\limits_{k \rightarrow\infty} \frac{1}{\varepsilon_k}
{\mathbb E}_{x_k}
 \bigl [   \bigl(1 - e^{-L (\newstop^{\varepsilon_k}_1)}   \bigr)
\mathbb{I}_{\{\tau
^0 <
\tau^1\}}   \bigr] \nonumber\\[-8pt]\\[-8pt]
  && \qquad = \operatorname{\lim\inf}\limits_{k \rightarrow\infty}
\frac
{1}{\varepsilon_k}
{\mathbb E}
  \bigl[  \bigl (1 - e^{-\newlk  (\newstop^{\varepsilon_k}_{(k),1}
)}   \bigr) \mathbb{I}_{\{\tau^{0}_{(k)} <
\tau^{1}_{(k)}\}}   \bigr],\nonumber
\end{eqnarray}
where $\tau_{(k)}^\varepsilon$ and $\tau^{k, \varepsilon}$ are
defined as follows:
%
\begin{eqnarray}
\label{def-tauk}
\tau_{(k)}^\varepsilon
& \doteq& \inf \bigl \{t \geq0 \dvtx \newzk(t)
\in H_\varepsilon
  \bigr\}, \nonumber\\[-8pt]\\[-8pt]
\tau^{k, \varepsilon}
& \doteq& \inf \bigl \{ t \geq0 \dvtx Z^k (t) \in
H_{\varepsilon/\varepsilon_k}
  \bigr\},\nonumber
\end{eqnarray}
and recall the definition of $\newstop_{(k),1}^{\varepsilon}$ given in
(\ref{def-newstopk1}).
Assume, without loss\vspace*{1pt} of generality, that $k$ is large enough so that
$\varepsilon_k < 1$.
Then, for each $x \geq0$, applying the mean value theorem to the function
$f_x(\varepsilon) = 1 - e^{-\varepsilon x}$, we infer that for $x \geq0$,
there exists $\varepsilon_k^* = \varepsilon_k^*(x) \in(0,\varepsilon
_k)$ such that
\[
\frac{1 - e^{-\varepsilon_k x}}{\varepsilon_k} = x e^{-\varepsilon
_k^* x}
\geq
x e^{-x}.
\]
Using the above inequality along with the equalities
$L_{(k)}(\newstop_{(k),1}^{\varepsilon_k}) = \varepsilon_k L^k
(\newtime_{k, 1})$,
$\varepsilon_k^2\tau^{k,0} = \tau_{(k)}^0$ and $\varepsilon_k^2\tau
^{k,1} = \tau
_{(k)}^1$, which hold due to
the scaling\vspace*{1pt} relations (\ref{def-scaling}) and (\ref{nscaling}), we have
for all $k$ sufficiently large,
\begin{eqnarray*}
\frac{1}{\varepsilon_k} {\mathbb E}  \bigl[  \bigl(1-e^{-\newlk
(\newstop
_{(k),1}^{\varepsilon_k})}  \bigr)
\mathbb{I}_{\{\tau^{0}_{(k)} < \tau^{1}_{(k)}\}}  \bigr]
&
= & {\mathbb E} \biggl [  \biggl(
\frac{1-e^{-\varepsilon_k L^k (\newtime_{k,1})}}{\varepsilon
_k}  \biggr)
\mathbb{I}_{\{\tau^{k,0}<\tau^{k,1}\}}   \biggr]\\
\\
& \geq&
{\mathbb E} \bigl [L^k(\newstop^{k,1}_1)e^{-L^k(\newstop^{k,1}_1)}
\mathbb{I}_{\{\tau^{k,0} < \tau^{k,1}\}}  \bigr].
\end{eqnarray*}
Comparing this with (\ref{eq-lbound2}) and (\ref{equivP}), it is
clear that to prove the lemma it
suffices to show that there exists $\tilde{C} > 0$ such that
%
\begin{equation}
\label{lem-toshow1}
\operatorname{\lim\inf}\limits_{k \rightarrow\infty} {\mathbb E} \bigl [L^k(\newtime
_{k,1})e^{-L^k(\newtime_{k,1})}
\mathbb{I}_{\{\tau^{k,0} < \tau^{k,1}\}}  \bigr] \geq\tilde{C}.
\end{equation}

Let $\overline{X} = x_* + \sigma(0)B$, where $B$ is standard Brownian motion,
and let $(\overline{Z}, \overline{Y})$ satisfy the ESP for
$\overline{X}$, as in Lemma \ref{lem-0}. Then,
since $x_k/\varepsilon_k \rightarrow x_* \in H_1$, by Lemma~\ref
{lem-0}(2) it follows
that there exist $r \in(1/2,1)$ and $R
\in(1, 2)$ such that as $k \rightarrow\infty$,
\[
\max_{i=1, \ldots, J} \sup_{s \in[0,\theta_{r,R}^k]} Y_i^k(s)
\Rightarrow
\max_{i=1, \ldots, J} \sup_{s \in[0,\overline{\theta}_{r,R}]}
\overline{Y}_i (s),
\]
where recall the definitions of
$\theta_{r,R}^k$ and $\overline{\theta}_{r,R}$ given in
(\ref{def-thetaeps}) and (\ref{def-thetalr}), respectively.
By the Portmanteau theorem,
this implies that
%
\begin{eqnarray*}
 {\mathbb P}  \Bigl( \max_{i=1,\ldots, J} \sup_{t \in[0,
\overline
{\theta}_{r, R}]} \overline{Y}_i(t) > \delta  \Bigr)
& \leq&
 \operatorname{\lim\inf}\limits_{k \rightarrow\infty} {\mathbb
P}  \Bigl( \max
_{i=1,\ldots,
J} \sup_{t \in[0, \theta^k_{r, R}]}
Y^k_i (t) > \delta  \Bigr) \\
& \leq&
 \operatorname{\lim\inf}\limits_{k \rightarrow\infty} {\mathbb
P}  \bigl( L^k
(\theta
^k_{r,R}  ) > \delta  \bigr).
\end{eqnarray*}
Together with the fact that property 3 of Corollary \ref{cor-0}
implies that
there exists $\delta> 0$ such that
\[
{\mathbb P}  \Bigl( \max_{i=1,\ldots, J} \sup_{t \in[0, \overline
{\theta
}_{r, R}]} \overline{Y}_i(t) > \delta  \Bigr)
> 2 \delta,
\]
and the inequality
$\newtime_{k,1} \geq\theta^{k}_{r, R}$ for all $k$,
it follows that there exists $K < \infty$ such that
%
\begin{equation}
\label{eq-one}
{\mathbb P} \bigl ( L^k   ( \newtime_{k,1}   ) > \delta
\bigr) \geq
\delta,\qquad
  k \geq K.
\end{equation}
Next, choose $r'\in(0,1/2)$ and $R' \in(2, \infty)$
and note that $\newtime_{k,1} \leq\theta^k_{r', R'}$
because $Z^k(0) \in H_1$ and $N_1(Z^k(0)) = H_{1/2} \cup H_2$.
Hence, property 1
of Corollary \ref{cor-0} implies that there exists $c < \infty$
such that
%
\begin{equation}
\label{eq-two}
\sup_{k \in{\mathbb N}} {\mathbb P}  \bigl(L^k (\theta^k_{r',
R'}) <c
\bigr) \geq
{\mathbb P} \Bigl ( \sup_{k\in{\mathbb N}} L^k(\theta^k_{r',
R'}) <c
\Bigr)\geq
1-\frac{\delta}{4}.
\end{equation}
On the other hand, since ${\mathbb P}_{Z(\tau^{\varepsilon_k})} (\tau
^0 \geq
\tau
^1) = \varepsilon_k$
by
Lemma \ref{lem-00} and $\varepsilon_k \rightarrow0$ as $k
\rightarrow\infty$, we have
\[
\lim_{k \rightarrow\infty} {\mathbb P}(\tau^{k,0}< \tau^{k,1})=
\lim_{k \rightarrow\infty} {\mathbb P}\bigl(\tau^0_{(k)}<\tau^1_{(k)}\bigr) =
\lim_{k \rightarrow\infty}   ( 1 - \varepsilon_k   ) = 1.
\]
Hence, by choosing $K < \infty$ larger if necessary, we can assume that
%
\begin{equation}
\label{eq-three}
{\mathbb P}(\tau^{k,0} < \tau^{k,1} ) \geq1 - \frac{\delta}{4},\qquad
  k
\geq K.
\end{equation}
Now, define the set
\[
S_k \doteq  \bigl\{\tau^{k,0} < \tau^{k,1}, e^{-L^k (\newtime
_{k,1})} \geq e^{-c},
L^{k}(\newtime_{k,1}) > \delta  \bigr\}.
\]
Then (\ref{eq-one}), (\ref{eq-two}) and (\ref{eq-three}), together
show that for $k \geq K$, ${\mathbb P}(S_k) \geq\delta/2$.
Therefore, for all $k \geq K$,
\begin{eqnarray*}
 &&{\mathbb E}  \bigl[L^k(\newtime_{k,1})e^{-L^k(\newtime_{k,1})}
\mathbb{I}_{\{\tau^{k,0} < \tau^{k,1}\}}  \bigr] \\
 &&\qquad \geq {\mathbb E}  \bigl[L^k(\newtime
_{k,1})e^{-L^k(\newtime_{k,1})}
\mathbb{I}_{\{\tau^{k,0} < \tau^{k,1}\}} \mathbb{I}_{S_k}   \bigr]
\geq\delta e^{-c} \frac{\delta}{2},
\end{eqnarray*}
and so
(\ref{lem-toshow1}) holds with $\tilde{C} = \delta^2 e^{-c}/2$.
This completes the proof of the lemma.
\end{pf}

\begin{lemma}
\label{lem-1}
Suppose $b \equiv0$.
For every $n \in{\mathbb N}$,
\[
\lim_{k \rightarrow\infty} \sup_{x \in H_{\varepsilon_k}} {\mathbb E}_x
  \bigl[ 1 - e^{-L(\newstop^{\varepsilon_k}_{n})}   \bigr] = 0.
\]
\end{lemma}
\begin{pf}
Fix $n \in{\mathbb N}$.
We prove the lemma using an argument by contradiction.
Suppose that there exists $\delta_0 > 0$ and a subsequence, which
we denote again by $\{\varepsilon_k\}_{k \in{\mathbb N}}$, such that
$\varepsilon_k \downarrow0$ as $k \rightarrow\infty$ and
for every $k \in{\mathbb N}$,
\[
\sup_{x\in H_{\varepsilon_k}}
{\mathbb E}_x \bigl [1-e^{-L(\newstop^{\varepsilon_k}_{n})}  \bigr]
\geq
\delta_0.
\]
{F}or each $k \in{\mathbb N}$, let $x_k \in H_{\varepsilon_k}$ be
such that
%
\begin{equation}
\label{contra1}
{\mathbb E}_{x_k}  \bigl[ 1-e^{-L(\newstop_{n}^{\varepsilon_k})}
\bigr] \geq
\frac{\delta_0}{2}.
\end{equation}
Since,
the law of $(\newzk,\newyk,\newlk)$ under ${\mathbb P}$ is the same
as the
law of $(Z,Y,L)$ under ${\mathbb P}_{x_k}$,
(\ref{contra1}) is equivalent to the inequality
%
\begin{equation}
\label{contra10}
{\mathbb E}  \bigl[ 1-e^{-\newlk(\newstop_{(k),n}^{\varepsilon
_k})}  \bigr]
\geq
\frac{\delta_0}{2}.
\end{equation}
%

The scaling relations in (\ref{def-scaling}) and (\ref{nscaling})
show that
%
\begin{equation}
\label{step-inter}
{\mathbb E} \bigl [ 1 -e^{-\newlk(\newstop_{(k),n}^{\varepsilon_k})}
  \bigr] =
{\mathbb E} \bigl [1 - e^{-\varepsilon_k L^k   (\newstop_{n}^{k,1}
  )}
  \bigr].
\end{equation}
Moreover, since
$Z^k(0) = x_k/\varepsilon_k \in H_1$, it follows that
$\langle Z^k(t), \vec{{\mathbf v}}\rangle\in  [ 2^{-n}, 2^n
  ]$ for $t \in  [0, \newtime_{k,n}  ]$,
Therefore, there exist $0<r<2^{-n}$ and $R> 2^n$ such that
$\newtime_{k,n} \leq\theta^k_{r, R}$,
where $\theta^k_{r, R}$ is defined in (\ref{def-thetaeps}).
As a result, we conclude that
\[
{\mathbb E}  \bigl[ 1-e^{-\varepsilon_k
L^k(\newtime^{k,n})}  \bigr]
\leq{\mathbb E} \bigl [1-e^{-\varepsilon_k
L^k(\theta^k_{r, R})}  \bigr] \rightarrow0 \qquad \mbox{as } k
\rightarrow\infty,
\]
where the last limit holds due to the weak convergence
$\varepsilon_k L^k (\theta^k_{r, R}) \Rightarrow0$ established
in Corollary \ref{cor-0},
and the fact that $x \mapsto1 - e^{-x}$ is a bounded
continuous function.
When combined with (\ref{step-inter}), this
contradicts (\ref{contra1}) and thus proves the lemma.
\end{pf}

We are now in a position to complete
the proof of Proposition \ref{prop-mainstep}.

\begin{pf*}{Proof of Proposition \ref{prop-mainstep}}
First, observe that by Lemma \ref{lem-2}, there exists $C > 0$ and $K
< \infty$ such that for all
$k \geq K$, the relation
\[
\inf_{x \in H_{\varepsilon_k}} {\mathbb E}_x   \bigl[
  \bigl(1-e^{-L(\newstop_1^{\varepsilon_k})}  \bigr) \mathbb{I}_{\{
\tau^0 < \tau
^1\}}   \bigr] \geq
\frac{C}{2}\varepsilon_k
\]
is satisfied.
Together with the fact that $Z(\tau^{\varepsilon_k})
\in H_{\varepsilon_k}$ and, for any $x \in H_{\varepsilon_k}$,
${\mathbb P}_x$-a.s.,
%
\begin{equation}
\label{eq-denom4}
\langle Z(\newstop^{\varepsilon_k}_{n-1}), \vec{{\mathbf v}}\rangle
\leq2^{n-1} \varepsilon_k,
\end{equation}
implies that for all $k$ large enough so that
$\varepsilon_k < 2^{-(n-1)} \varepsilon_0$ and for $m = 1, \ldots, n$,
%
\begin{eqnarray}\label{eq-denom5}
&&{\mathbb E}_{Z(\tau^{\varepsilon_k})}  \bigl [ e^{-L(\newstop
^{\varepsilon
_k}_{m-1})} {\mathbb E}_{Z(\newstop^{\varepsilon_k}_{m-1})}
\bigl[\bigl(1-e^{-L(\newstop
_1^{\varepsilon_k})}\bigr)\mathbb{I}_{\{\tau^0 <
\tau^1\}}  \bigr]   \bigr] \nonumber\\[-8pt]\\[-8pt]
 &&\qquad\geq\frac{C}{2} {\mathbb E}_{Z(\tau^{\varepsilon_k})}   \bigl[
e^{-L(\newstop^{\varepsilon_k}_{m-1})} \langle
Z(\newstop^{\varepsilon_k}_{m-1}), \vec{{\mathbf v}}\rangle  \bigr].\nonumber
\end{eqnarray}
When combined with (\ref{eq-denom2}) and (\ref{eq-denom3}), this
shows that
%
\begin{eqnarray}
\label{insert1}
&& {\mathbb E}_0   \bigl[ {\mathbb E}_{Z(\tau^{\varepsilon_k})}
  \bigl[ \bigl (1-e^{-L(\tau^0)}  \bigr) \mathbb{I}_{\{\tau^0< \tau
^1\}}
  \bigr]  \bigr]\nonumber\\[-8pt]\\[-8pt]
 &&\qquad \geq
\frac{C}{2}\sum_{m=1}^n {\mathbb E}_{Z(\tau^{\varepsilon_k})}
  \bigl[
e^{-L(\newstop^{\varepsilon_k}_{m-1})} \langle
Z(\newstop^{\varepsilon_k}_{m-1}), \vec{{\mathbf v}}\rangle  \bigr].\nonumber
\end{eqnarray}
Each summand on the right-hand side
can be rewritten in the more convenient form
\begin{eqnarray*}
&&{\mathbb E}_{Z(\tau^{\varepsilon_k})}   \bigl[ e^{-L(\newstop
^{\varepsilon
_k}_{m-1})} \langle
Z(\newstop^{\varepsilon_k}_{m-1}), \vec{{\mathbf v}}\rangle  \bigr]\\
&&\qquad =
{\mathbb E}_{Z(\tau^{\varepsilon_k})}   [ \langle
Z(\newstop^{\varepsilon_k}_{m-1}), \vec{{\mathbf v}}\rangle  ]
 - {\mathbb E}_{Z(\tau^{\varepsilon_k})}   \bigl[  \bigl( 1 -
e^{-L(\newstop
^{\varepsilon_k}_{m-1})}  \bigr) \langle
Z(\newstop^{\varepsilon_k}_{m-1}), \vec{{\mathbf v}}\rangle  \bigr].
\end{eqnarray*}
Since $b \equiv0$, Lemma \ref{lem-00} and the uniform bound (\ref{eq-denom4})
show that $\langle Z, \vec{{\mathbf v}}\rangle$ is a
martingale on $[0,\newstop_n^{\varepsilon_k}]$. In addition, because
$\newstop^{\varepsilon_k}_{m-1} \leq\newstop_n^{\varepsilon_k}$
and $\langle Z(\tau^{\varepsilon_k}), \vec{{\mathbf v}}\rangle=
\varepsilon_k$, it follows that
\[
{\mathbb E}_0 \bigl [{\mathbb E}_{Z(\tau^{\varepsilon_k})}   [
\langle
Z(\newstop^{\varepsilon_k}_{m-1}), \vec{{\mathbf v}}\rangle  ]
  \bigr] =
{\mathbb E}_0
[\varepsilon_k] = \varepsilon_k.
\]
Furthermore, by (\ref{eq-denom4}), Lemma \ref{lem-1} and the bounded
convergence theorem, we have
for any $n \in{\mathbb N}$ and $m = 1, \ldots, n$,
\begin{eqnarray*}
 &&\operatorname{\lim\sup}\limits_{k \rightarrow\infty} \frac{1}{\varepsilon_k}
{\mathbb E}_0
  \bigl[{\mathbb E}_{Z(\tau^{\varepsilon_k})}  \bigl[
\bigl(1-e^{-L(\newstop
^{\varepsilon_k}_{m-1})}  \bigr) \langle
Z(\newstop^{\varepsilon_k}_{m-1}), \vec{{\mathbf v}}\rangle
\bigr]  \bigr]\\
&& \qquad
 \leq 2^{n-1}\lim_{k \rightarrow\infty} {\mathbb E}_0
  \Bigl[
\sup_{x\in H_{\varepsilon_k}}
{\mathbb E}_x  \bigl[
1-e^{-L(\newstop^{\varepsilon_k}_{m-1})}   \bigr]   \Bigr]  =  0.
\end{eqnarray*}
Combining the last three assertions, we see that for every $n \in
{\mathbb N}$
and $m = 1, \ldots, n$,
\[
\operatorname{\lim\inf}\limits_{k \rightarrow\infty} \frac{1}{\varepsilon_k} {\mathbb
E}_0  \bigl[
{\mathbb E}_{Z(\tau
^{\varepsilon_k})}   \bigl[ e^{-L(\newstop^{\varepsilon_k}_{m-1})}
\langle
Z(\newstop^{\varepsilon_k}_{m-1}), \vec{{\mathbf v}}\rangle  \bigr]
  \bigr] = 1.
\]
Together
with (\ref{insert1}), this shows that for every $n \in{\mathbb N}$,
\[
\operatorname{\lim\inf}\limits_{k \rightarrow\infty} \frac{1}{\varepsilon_k} {\mathbb
E}_0   \bigl[
{\mathbb E}_{Z(\tau
^{\varepsilon_k})}
  \bigl[ \bigl (1-e^{-L(\tau^0)}  \bigr) \mathbb{I}_{\{\tau^0< \tau
^1\}}
  \bigr]  \bigr] \geq\frac{nC}{2}.
\]
Taking the limit as $n \rightarrow\infty$, we obtain (\ref{toshow-denom}),
thus completing
the proof of the proposition.
\end{pf*}

\subsection{The general drift case}
\label{subs-gendrift}

In this section, we establish Theorem \ref{th-mainstep}.
Specifically, we use a Girsanov transformation to
generalize the case of zero drift,
established in Proposition \ref{prop-mainstep},
to arbitrary Lipschitz drifts with linear growth,
as specified in property 3 of Definition \ref{def-gpssder}.
As before, let $Z$ be the unique strong solution
to the {Class ${\mathcaligr A}$} SDER, which exists by
Theorem \ref{th-SDER},
and let $\tau^1$ be the first
hitting time to $H_1$, as defined in (\ref{def-tauve}).
We begin with a simple lemma that shows that $\tau^1$ is finite with
positive ${\mathbb P}_0$ probability.

\begin{lemma}
\label{lem-finite} We have
%
\begin{equation}
\label{toshow-drift2}
{\mathbb P}_0   ( \tau^1 < \infty  ) > 0.
\end{equation}
Moreover, if $\inf_{x: \langle x, \vec{{\mathbf v}}\rangle\leq1}
\langle b(x),
\vec{{\mathbf v}}\rangle\geq0$, then
%
\begin{equation}
\label{toshow-drift}
{\mathbb P}_0  (\tau^1 < \infty  ) =1.
\end{equation}
\end{lemma}
\begin{pf}
Recall the definition of $X$ and $M$ given in (\ref{def-x}) and
(\ref{def-MA1}).
By Theorem \ref{th-SDER}, we know that ${\mathbb P}$-a.s., $Z$ satisfies
the ESP for $X$. Hence, by Lemma \ref{lem-vset} it follows that
$\widehat{Z} = \Gamma_1 (\widehat{X})$, where $\Gamma_1$ is the
$1$-dimensional
Skorokhod map and, for $H = Z, M, X$, we define
$\widehat H \doteq\langle H, \vec{{\mathbf v}}\rangle$.
Let $T(t) \doteq\inf\{ s \geq0 \dvtx \langle\widehat{M} \rangle_s >
t\}$.
Then, due to the uniform ellipticity of $a$, $T$ is strictly increasing and,
since $\widehat{M}$ is a continuous martingale,
$\widehat{M}(T(\cdot))$ is a $1$-dimensional
Brownian motion. In turn, this implies $\widehat{Z}$ is
a one-dimensional reflected Brownian motion with drift
\[
\int^t_0\langle b(Z(T(s))), \vec{{\mathbf v}}\rangle\, dT(s) = \int^t_0
\langle
b(Z(T(s))), \vec{{\mathbf v}}\rangle
\frac{1}{\vec{{\mathbf v}}^T a(Z(s))\vec{{\mathbf v}}} \, ds.
\]
Since $\langle b(x), \vec{{\mathbf v}}\rangle/\vec{{\mathbf v}}^T
a(x)\vec{{\mathbf v}}$ is
continuous on $G$, there exists $\kappa\in(-\infty, \infty)$ such that
\[
\frac{\langle b(x), \vec{{\mathbf v}}\rangle}{\vec{{\mathbf v}}^T
a(x)\vec{{\mathbf v}}} >
\kappa\qquad \mbox{for all }
x \in G, \langle x, \vec{{\mathbf v}}\rangle\leq1.
\]
Consider the process $\tilde{X}$ defined by $\tilde{X}(t) \doteq
\kappa t +
M(T(t))$ for $t \in[0,\infty)$,
and let $\tilde{Z} \doteq\Gamma_1(\tilde{X})$
be a one-dimensional reflected Brownian motion with constant drift
$\kappa$.
Then $\widehat{X}(T(t)) - \widehat{X}(T(s)) \geq\tilde{X}(t) -
\tilde{X}(s)$ for every $0 \leq s \leq t$, and so
the comparison principle for $\Gamma_1$ (see, e.g., equation
(4.1) in Lemma 4.1 of \cite{krulehramshr1}) shows that $\widehat
{Z}(T(t)) \geq\tilde{Z}(t)$ for every $t \in[0,\widehat{\tau}^1]$, where
\[
\widehat{\tau}^1 \doteq\inf\{ t > 0 \dvtx \widehat{Z} (T(t)) = 1 \}.
\]
Since $T(\widehat{\tau}^1)=\tau^1$, it follows that
\[
{\mathbb P}_0 \bigl(\widehat{Z} \bigl(T(t) \wedge\tau^1\bigr)\geq\tilde{Z}
(t\wedge
\widehat{\tau}^1) \mbox{ for all } t \geq0 \bigr) = 1.
\]
Since $T$ is strictly increasing, we have $\tau^1 = \infty$ if and
only if
$\widehat{\tau}^1 = \infty$. Therefore, on the set $\{\tau^1 =
\infty\}$, we must have
\[
\tilde{Z}(t) \leq\widehat{Z}(T(t)) < 1 \qquad \mbox{for all } t
\in[0,\infty).
\]
However, $\tilde{Z}$ will hit $1$ with positive ${\mathbb P}_0$ probability,
and in fact will hit $1$
${\mathbb P}_0$-a.s. if $\kappa\geq0$ (see, e.g., page 197 of
\cite{karshrbook}), and so the same is true of $\hat{Z}(T(\cdot))$.
This implies both (\ref{toshow-drift2}) and (\ref{toshow-drift}),
and so the proof of the lemma is complete.
\end{pf}

\begin{pf*}{Proof of Theorem \ref{th-mainstep}}
The uniform ellipticity of $a(\cdot)$ ensures that
$a^{-1}(\cdot)$ exists. Let $\mu\doteq
-\sigma^T a^{-1} b$, note that $\mu^T \mu= b^T a b$, and define
%
\begin{equation}\label{D}\qquad
   D(t) \doteq\exp  \biggl\{\int^t_0 \mu(Z(s)) \, dB(s) -
\frac{1}{2} \int^t_0 b^T (Z(s)) a(Z(s)) b (Z(s)) \, ds   \biggr\}
\end{equation}
for $t \in[0,\infty)$.
Property 3 of Definition \ref{def-gpssder}
guarantees that $\mu$ has at most linear growth,
and so, as is well-known, $\{D(t), {\mathcaligr F}_t\}$ is a
martingale (see,
e.g., Corollary~5.16 of \cite{karshrbook}).

Fix $T < \infty$. Define a new probability measure ${\mathbb Q}_0$ on
$(\Omega
, {\mathcaligr F}, \{{\mathcaligr F}_T\})$ by setting
\[
{\mathbb Q}_0 (A) = {\mathbb E}  [ D(T) \mathbb{I}_{A}  ]
\qquad\mbox{for } A
\in{\mathcaligr F}_{T}.
\]
Define
\[
\tilde{B}(t) \doteq B(t) + \int^t_0 \sigma^T
(Z(s))a^{-1}(Z(s))b(s) \, ds,\qquad   t \in[0,T].
\]
By Girsanov's theorem (see Theorem 5.1 of \cite{karshrbook}), under
${\mathbb Q}_0$,
$\{\tilde{B}_t, {\mathcaligr F}_t\}_{t \in[0,T]}$ is a Brownian
motion and
\[
Z(t) = \int^{t}_0 \sigma(Z(s)) \, d \tilde{B} (s) + Y(t),\qquad   t
\in[0,T],
\]
where $(Z,Y)$ satisfy the ESP pathwise for $Z-Y$.
Since, under ${\mathbb Q}_0$, $Z$ is the solution to a {Class
${\mathcaligr A}$} SDER
with no drift,
by Proposition \ref{prop-mainstep}, it follows that
\[
{\mathbb Q}_0  \bigl( L(\tau^1)<\infty, \tau^1 \leq T   \bigr) = 0.
\]
Since ${\mathbb P}_0 \ll{\mathbb Q}_0$ [with $d{\mathbb P}_0/d{\mathbb
Q}_0$ $= D^{-1} (T)$ on
${\mathcaligr F}_{T}$],
this implies
\[
{\mathbb P}_0   \bigl( L(\tau^1) <\infty, \tau^1 \leq T   \bigr) = 0.
\]
Since $T < \infty$ is arbitrary,
sending $T \rightarrow\infty$ (along a countable sequence), we
conclude that
\[
{\mathbb P}_0   \bigl( L(\tau^1) < \infty,  \tau^1 < \infty  \bigr)
= 0.
\]
However, ${\mathbb P}_0 (\tau^1 < \infty) > 0$ by Lemma \ref
{lem-finite}. Hence,
${\mathbb P}_0(L(\tau^1) = \infty, \tau^1 < \infty) > 0$, which in turn
implies that there exists
$T < \infty$ such that ${\mathbb P}_0(L(T) = \infty) > 0$, which proves
Theorem \ref{th-mainstep}.
In addition, note that if $\inf_{x\in G:\langle x, \vec{{\mathbf
v}}\rangle
\leq1}
\langle b(x), \vec{{\mathbf v}}\rangle\geq0$, then ${\mathbb P}_0(
\tau^1 <
\infty)=1$
and so we in fact
have ${\mathbb P}_0 (L(\tau^1) =\infty) = 1$.
\end{pf*}

\subsection{The semimartingale property for $Z$}
\label{subs-main1pf}

Recall from Theorem \ref{th-SDER} that the process $Z$ has the decomposition
$Z = M + A$, where
%
\begin{equation}\label{MA} M =\int_0^\cdot
\sigma(Z(s)) \, d B(s),\qquad   A = \int_0^\cdot b(Z(s)) \, ds + Y,
\end{equation}
and $Y$ is the constraining term associated with the ESP.
$M$ is clearly a (local) martingale, but Theorem \ref{th-mainstep} shows
that $Y$ is not ${\mathbb P}$-a.s. of finite variation on bounded intervals.
However,
as mentioned earlier, Theorem \ref{th-mainstep} does not immediately
imply that $Z$ is not a semimartingale because we do not know a priori
that the above decomposition must be the Doob decomposition of $Z$ if
it were
a semimartingale.
In Proposition \ref{prop-smg} below, we show that the latter
statement is indeed true, thus showing that $Z$ is not a semimartingale.

\begin{prop}
\label{prop-smg}
If $Z$ were a semimartingale, then its Doob decomposition must be $Z =
M + A$.
\end{prop}
\begin{pf}
Suppose that $Z$ is a semimartingale, and let its (unique) Doob decomposition
take the form
\[
Z = \tilde{M} +\tilde{A},
\]
where $\tilde{M}$ is an $\{\mathcaligr F_t\}$-adapted continuous local
martingale and $\tilde{A}$ is an $\{\mathcaligr F_t\}$-adapted continuous,
process with ${\mathbb P}$-a.s. finite variation on bounded intervals.

Fix $R < \infty$ and let $\theta_R\doteq\inf\{t \geq0 \dvtx |M(t)|
\geq R\}$.
For each $\varepsilon>0$,
define two sequences of stopping times $\{\tau^\varepsilon_n\}_{n \in
{\mathbb N}}$ and
$\{\xi^\varepsilon_n\}_{n \in{\mathbb N}}$ as follows: $\xi
^\varepsilon_0 \doteq
0$ and for
$n \in{\mathbb N}$, let
\begin{eqnarray*}
\tau^{\varepsilon}_n & \doteq& \inf  \{ t \geq\xi
_{n-1}^\varepsilon \dvtx Z(t)
\in H_\varepsilon  \} \wedge\theta_R, \\
\xi^{\varepsilon}_n & \doteq& \inf  \{ t \geq\tau
_n^\varepsilon \dvtx
Z(t) \in H_{\varepsilon/2}   \} \wedge\theta_R.
\end{eqnarray*}
(For notational conciseness, we have suppressed the dependence of these
stopping times on $R$.)
By uniqueness of the Doob decomposition, clearly
$Z(\cdot\wedge\xi^\varepsilon_n) - Z(\cdot\wedge
\tau^{\varepsilon}_n)$ is an $\{\mathcaligr F_t\}$-adapted
semimartingale, with
Doob decomposition
\[
Z(t\wedge\xi^\varepsilon_n) - Z(t\wedge
\tau^{\varepsilon}_n)
= \tilde M(t \wedge
\xi^\varepsilon_n) -\tilde M(t \wedge\tau^{\varepsilon}_n) +
\tilde A(t\wedge\xi^\varepsilon_n) - \tilde A(t\wedge
\tau^{\varepsilon}_n)
\]
On the other hand, due to the identity
$Z = M + A = \tilde{M} + \tilde{A}$, we also have
\[
Z(t\wedge\xi^\varepsilon_n) - Z(t\wedge
\tau^{\varepsilon}_n) =  M(t \wedge
\xi^\varepsilon_n) -M(t \wedge\tau^{\varepsilon}_n) + A(t\wedge
\xi^\varepsilon
_n) -
A(t\wedge\tau^{\varepsilon}_n).
\]
Since $M$ is an $\{\mathcaligr F_t\}$-adapted continuous (local)
martingale, and
$M$ is uniformly bounded on $[0,\theta_R]$,
the stopped processes $M(\cdot\wedge\xi^\varepsilon_n)$ and
$M(\cdot\wedge
\tau^{\varepsilon}_n)$ are
$\{\mathcaligr F_t\}$-adapted continuous martingales.
Hence, $M(\cdot\wedge\xi^\varepsilon_n)-M(\cdot\wedge\tau
^{\varepsilon
}_n)$ is also an $\{\mathcaligr
F_t\}$-adapted continuous martingale.
Moreover, Theorem \ref{th-SDER} implies that $Y(\cdot\wedge\xi
^\varepsilon
_n) - Y(\cdot\wedge
\tau^{\varepsilon}_n)$ has ${\mathbb P}$-a.s. finite variation on each
bounded time
interval. Since $A = Y + \int_0^\cdot b(Z(s)) \, ds$,
$A(\cdot\wedge\xi^\varepsilon_n) - A(\cdot\wedge\tau^\varepsilon
_n)$ also
has ${\mathbb P}$-a.s. finite variation on each bounded time interval.
By uniqueness of the Doob
decomposition, we conclude that for every $\varepsilon> 0$ and $t\in[0,
\infty)$,
\[
M(t\wedge\xi^\varepsilon_n) - M(t \wedge
\tau^{\varepsilon}_n) = \tilde{M}(t\wedge\xi^\varepsilon_n)
- \tilde{M}(t\wedge\tau^{\varepsilon}_n).
\]
Summing over $n \in{\mathbb N}$ on both sides of the last equation, we obtain
%
\begin{equation}
\label{equality-1}
\sum^{\infty}_{n=1} \bigl(M(t\wedge
\xi^\varepsilon_n) - M(t\wedge\tau^{\varepsilon}_n)\bigr) =
\sum^{\infty}_{n=1} \bigl(\tilde{M}(t
\wedge\xi^\varepsilon_n)- \tilde{M}(t \wedge
\tau^{\varepsilon}_n)\bigr).
\end{equation}

On the other hand, because ${\mathbb P}$-a.s., $M(0) = 0$ and $\xi
_n^\varepsilon\rightarrow
\theta_R$ as $n \rightarrow\infty$,
we can write $M(t\wedge\theta_R)$ as a telescopic sum:
\[
M(t \wedge\theta_{R}) = \sum^{\infty}_{n=1}\bigl(M(t\wedge
\xi^\varepsilon_n) - M(t\wedge\xi^{\varepsilon}_{n-1})\bigr), \qquad  t
\in[0,\infty).
\]
Next, observe that
\begin{eqnarray*}
&&M(t \wedge\theta_{R}) - \sum^{\infty}_{n=1}\bigl(M(t\wedge
\xi^\varepsilon_n) - M(t\wedge\tau^{\varepsilon}_n)\bigr)\\
&&\qquad= \sum^{\infty}_{n=1} \bigl(M(t
\wedge\tau^{\varepsilon}_n) - M(t \wedge
\xi^\varepsilon_{n-1})\bigr) \\
 &&\qquad =  \int^t_0
\sum^{\infty}_{n=1}\mathbb{I}_{(\xi^\varepsilon_{n-1}, \tau^{\varepsilon}_n]}
(s) \, dM(s) \\
 &&\qquad =  \int^t_0 \sum^{\infty}_{n=1}\mathbb{I}_{(\xi^\varepsilon_{n-1}, \tau^{\varepsilon}_n]}
(s) \mathbb{I}_{[0,\varepsilon]} (\langle\vec{{\mathbf v}},
Z(s)\rangle) \, dM (s),
\end{eqnarray*}
where the last equality holds because $\langle\vec{{\mathbf v}},
Z(s)\rangle
\leq\varepsilon$ for
$s\in(\xi^\varepsilon_{n-1}, \tau^{\varepsilon}_n]$.
When combined with Doob's maximal martingale inequality, this yields
\begin{eqnarray*}
&& {\mathbb E}  \Biggl[ \sup_{s\in[0, t]}   \Bigg| M(s
\wedge\theta
_{R} ) -
\sum^{\infty}_{n=1}\bigl(M(s\wedge\xi^\varepsilon_n) - M(s\wedge
\tau^{\varepsilon}_n)\bigr)  \Bigg |^2   \Biggr] \\
 && \qquad \leq4{\mathbb E} \Biggl [   \Bigg| M(t \wedge
\theta
_R)
- \sum^{\infty}_{n=1}\bigl(M(t\wedge\xi^\varepsilon_n) - M(t\wedge
\tau^{\varepsilon}_n)\bigr)  \Bigg |^2   \Biggr]\\
 && \qquad  = 4{\mathbb E}  \Biggl[  \Bigg| \int^t_0
\sum^{\infty}_{n=1}\mathbb{I}_{(\xi^\varepsilon_{n-1}, \tau^{\varepsilon
}_n]}(s) \mathbb{I}
_{[0,\varepsilon]} (\langle\vec{{\mathbf v}}, Z(s)\rangle)
\, dM(s)  \Bigg|^2  \Biggr]\\
&& \qquad\leq4{\mathbb E}  \biggl[ \int^t_0
\mathbb{I}_{[0,\varepsilon]} (\langle\vec{{\mathbf v}}, Z(s)\rangle
)   |a(Z(s))
| \, ds   \biggr].
\end{eqnarray*}
By Assumption \ref{as-locbd}, $a$ is bounded on the set $\{x\dvtx \langle
\vec{{\mathbf v}},
x\rangle\leq\varepsilon\}$. Hence, an application of
the bounded convergence theorem shows that
\[
 \lim_{\varepsilon\rightarrow0} {\mathbb E}  \biggl[ \int^t_0
\mathbb{I}_{\{\langle\vec{{\mathbf v}}, Z(s)\rangle
\leq\varepsilon\}}   |a(Z(s))  | \, ds   \biggr]
=  |a(0)| {\mathbb E} \biggl [ \int^t_0 \mathbb{I}_{\{
\langle\vec{{\mathbf v}},
Z(s)\rangle=0\}} \, ds
  \biggr]
= 0,
\]
where the last equality is a consequence of the fact that
$\langle\vec{{\mathbf v}}, Z\rangle$ is a uniformly elliptic one-dimensional
reflected diffusion (see Lemma \ref{lem-vset}) and consequently
spends zero Lebesgue time at the origin (see, e.g., page 90 of
\cite{freibook}).

An exactly analogous argument, with $\tilde{\theta}_{R} \doteq
\inf\{t \geq0 \dvtx |\tilde{M}|(t) \geq R\}$ and
$\tilde{\xi}_n^\varepsilon, \tilde{\tau}_n^\varepsilon$ defined
in a fashion analogous
to $\xi_n^\varepsilon, \tau_n^\varepsilon$, but with $\theta_R$
replaced by
$\tilde\theta_R$, shows that
\begin{eqnarray*}
 &&\lim_{\varepsilon\rightarrow0} {\mathbb E}  \Biggl[
\sup_{s\in[0, t]}
\Bigg|\tilde{M}
(s \wedge\tilde{\theta}_{R}) -
\sum^{\infty}_{n=1}\bigl(\tilde{M}(s\wedge\tilde{\xi}^\varepsilon_n)
- \tilde{M}(s\wedge\tilde{\tau}^{\varepsilon}_n)\bigr)   \Bigg|^2
  \Biggr] \\
&&\qquad \leq \lim_{\varepsilon\rightarrow0} 4 J^2 \sum_{i=1}^J
{\mathbb E}  \biggl[ \int^t_0 \mathbb{I}_{[0,\varepsilon]} (\langle
\vec{{\mathbf v}}, Z(s)
\rangle
) \, d\langle\tilde M_i \rangle(s)   \biggr] \\
&&\qquad =  4 J^2 \sum_{i=1}^J
{\mathbb E}  \biggl[\int^t_0 \mathbb{I}_{\{0\}} (\langle\vec{{\mathbf
v}}, Z(s)
\rangle) \,
d\langle\tilde M_i \rangle(s)   \biggr] \\
&&\qquad = 4 J^2 \sum_{i=1}^J
{\mathbb E} \biggl [\int^t_0 \mathbb{I}_{\{0\}} (Z_i(s)) \, d\langle
\tilde M_i
\rangle(s)
  \biggr],
\end{eqnarray*}
where
the last equality uses the property that
$Z_i(s)=0$ for every $i = 1, \ldots, J$ if and only if
$\langle\vec{{\mathbf v}}, Z(s) \rangle= 0$
(see property 2 of Definition \ref{def-gpssder}).
Due to the assumption that $\tilde{Z}_i$ is a semimartingale
with decomposition $\tilde{M}_i + \tilde{A}_i$,
the occupation times formula for continuous semimartingales (see, e.g.,
Corollary 1.6 in Chapter VI of \cite{RY}) and
the fact that the set $\{x\dvtx x_i =0\}$ has zero Lebesgue measure, we have,
${\mathbb P}$-a.s., for $i = 1, \ldots, J$,
\[
\int^t_0 \mathbb{I}_{\{0\}} (Z_i (s)) \, d\langle\tilde{M}_i
\rangle(s) =
\int^t_0 \mathbb{I}_{\{0\}} (Z_i (s)) \, d\langle Z_i \rangle(s) = 0.
\]
Combining the last four displays with (\ref{equality-1}), we conclude that
$M(t \wedge\theta_R) = \tilde{M} (t \wedge\tilde\theta
_{R} )$, ${\mathbb P}_0$-a.s., for every
$t \geq0$. This in turn implies that $\theta_R= \tilde\theta
_R$ ${\mathbb P}_0$-a.s.\vspace*{1pt}
Sending $R \rightarrow\infty$ and invoking
the continuity of both $M$ and $\tilde M$, we conclude that
$M = \tilde{M}$ ${\mathbb P}_0$-a.s.
In turn, this implies $A = \tilde{A}$, thus completing the proof of
the theorem.
\end{pf}

The proof of Theorem \ref{th-main1} is now a simple consequence of
Theorem \ref{th-mainstep}
and Proposition \ref{prop-smg}.

\begin{pf*}{Proof of Theorem \protect\ref{th-main1}}
If $Z$ were a semimartingale under ${\mathbb P}_0$,
then by Proposition \ref{prop-smg}, $Z = M + A$ is the Doob
decomposition for $Z$. In particular, this implies that
${\mathbb P}_0 (L(T) < \infty) = 1$ for every $T \in[0,\infty)$, where
recall that $L(T) = {\mathrm{Var}}_{[0,T]} Y$.
However, this contradicts the assertion of
Theorem \ref{th-mainstep} that there exists $T < \infty$ such that
${\mathbb P}_0 (L(T) = \infty) > 0$. Thus, we conclude that $Z$ is not a
semimartingale.
\end{pf*}

\begin{remark}
It is natural to expect
that similar, but somewhat more involved, arguments
could be used to show that the semimartingale property fails to
hold for a more general class of reflected diffusions in the
nonnegative orthant, in particular those
that arise as approximations of generalized processor
sharing networks (rather than just a single station, as considered in
\cite{ramrei1,ramrei2}).
Such diffusions would
satisfy properties 1, 2 and 4 of Definition \ref{def-gpssder}
but would have more complicated ${\mathcaligr V}$-sets
(see \cite{dupram00} for a description of the ESP associated
with such a network).
This is a subject for future work.
\end{remark}

\section{Dirichlet process characterization}
\label{sec-dirichlet}

This section is devoted to the proof of Theorem \ref{th-main2}.
Specifically, here we only assume that $(G,d(\cdot))$, $b(\cdot)$ and
$\sigma(\cdot)$ satisfy
Assumptions \ref{as-weak} and \ref{as-locbd}, and let
$(Z_t, B_t), (\Omega, {\mathcaligr F}, {\mathbb P}), \{{\mathcaligr
F}_t\}$ be
a Markov process that is a weak
solution to the associated SDER that satisfies Assumption
\ref{as-weakcont} for some constants $p > 1, q \geq2$ and
$K_T < \infty, T \in(0,\infty)$.
As usual, let $X$ be as defined in~(\ref{def-x}), and let
$Y = Z - X$, so that we can write
\[
Z(t)= Z(0) + \int_0^t b(Z(s))\, ds+ \int_0^t \sigma(Z(s)) \, d
B(s) +
Y(t), \qquad  t \in[0,\infty).
\]
Note that $\int_0^\cdot b(Z(s)) \, ds$ is a process of bounded
variation, and therefore
of bounded $p$-variation for any $p > 1$ by Remark \ref{zeroEn}.
As a result, in order to establish Theorem~\ref{th-main2}, it suffices
to show that
under ${\mathbb P}$-a.s., $Y$ has zero $p$-variation.

In Section \ref{subs-equiv}, we first show that it suffices to establish
a localized version (\ref{YzeroEa}) of the zero $p$-variation
condition on $Y$.
This is then used to prove Theorem~\ref{th-main2} in Section \ref
{subs-dirzero}.

\subsection{Localization}
\label{subs-equiv}

Fix $T>0$, let $\{\pi^n,n\geq1\}$ be a
sequence of partitions of $[0,T]$ such that $\Delta(\pi^n)\rightarrow0$
as $n\rightarrow\infty$. As mentioned above, to prove Theorem~\ref
{th-main2} we need to
establish the following result:
%
\begin{equation}
\label{YzeroE}
\sum_{t_i\in\pi^n}|Y(t_i)-Y(t_{i-1})|^p \stackrel{({\mathbb
P})}{\rightarrow} 0\qquad
\mbox{as } \Delta(\pi^n)\rightarrow\infty.
\end{equation}
{F}or each $m \in(0, \infty)$, let
%
\begin{equation}\label{zetam} \zeta^m \doteq\inf\{ t > 0 \dvtx   |
Z(t)
| \geq m\}.
\end{equation}
It is easy to see that ${\mathbb P}$-a.s., $\zeta^m\rightarrow\infty
$ as $m
\rightarrow\infty$.
We now show that the localized version, (\ref{YzeroEa}) below, is
equivalent to (\ref{YzeroE}).

\begin{lemma} \label{lem:Ytau}
The result (\ref{YzeroE}) holds if and only if for each $m \in
(0,\infty)$,
%
\begin{equation}
\label{YzeroEa}
\sum_{t_i\in\pi^n}|Y(t_i\wedge\zeta^m)-Y(t_{i-1}\wedge\zeta
^m)|^p \stackrel{({\mathbb P})}{\rightarrow} 0
\qquad\mbox{as } \Delta(\pi^n)\rightarrow0.
\end{equation}
\end{lemma}
\begin{pf} First, assume (\ref{YzeroEa}) holds for every $m \in(0,\infty)$.
Then, for every $m \in(0,\infty)$ and $\delta>0$,
\begin{eqnarray*}
 &&{\mathbb P}  \biggl(\sum_{t_i\in\pi
^n}|Y(t_i)-Y(t_{i-1})|^p \geq
\delta
  \biggr) \\
 &&\qquad\leq {\mathbb P}  \biggl(\sum_{t_i\in\pi
^n}|Y(t_i)-Y(t_{i-1})|^p
\geq\delta,\zeta^m>T   \biggr) + {\mathbb P}(\zeta^m\leq T) \\
&&\qquad =  {\mathbb P}  \biggl(\sum_{t_i\in\pi^n}|Y(t_i
\wedge\zeta
^m)-Y(t_{i-1} \wedge\zeta^m)|^p \geq\delta, \zeta^m>T   \biggr) +
{\mathbb P}(\zeta^m\leq T) \\
&&\qquad\leq {\mathbb
P}  \biggl(\sum
_{t_i\in\pi
^n}|Y(t_i\wedge\zeta^m)-Y(t_{i-1}\wedge\zeta^m)|^p \geq\delta
  \biggr) + {\mathbb P}(\zeta^m\leq T).
\end{eqnarray*}
Taking limits as $\Delta(\pi^n) \rightarrow0$, the first term on the
right-hand side vanishes due to~(\ref{YzeroEa}). Next,
sending $m \rightarrow\infty$, and using the fact that $\zeta
^m\rightarrow
\infty$
${\mathbb P}$-a.s., the second
term also vanishes, and so
we obtain (\ref{YzeroE}). This proves the ``if'' part of the result.

In order to prove the converse result,
suppose (\ref{YzeroE}) holds. Let $\theta^m_n \doteq\sup\{t_i \in
\pi^n \dvtx t_i\leq\zeta^m\}$, where $\theta_n^m \doteq T$ if the latter set is empty.
Then
\begin{eqnarray*}
 &&\sum_{t_i\in\pi^n}|Y(t_i\wedge\zeta
^m)-Y(t_{i-1}\wedge\zeta^m)|^p \\
  &&\qquad\leq\sum_{t_i\in\pi^n}|Y(t_i)-Y(t_{i-1})|^p
+|Y(\zeta^m \wedge T)-Y(\theta^m_n)|^p.
\end{eqnarray*}
Taking limits as $\Delta(\pi^n) \rightarrow0$, ${\mathbb P}$-a.s
the
last term vanishes since $|\zeta^m \wedge T -\theta^m_n| \leq\Delta
(\pi^n)$ and $Y$ is continuous.
Therefore, (\ref{YzeroEa}) follows from (\ref{YzeroE}).
\end{pf}

\subsection{The decomposition result}
\label{subs-dirzero}

For each $\varepsilon>0$, recursively
define two sequences of stopping times $\{\tau^\varepsilon_n\}_{n \in
{\mathbb N}}$ and
$\{\xi^\varepsilon_n\}_{n \in{\mathbb N}}$ as follows: $\xi
^\varepsilon_0 \doteq
0$ and for
$n \in{\mathbb N}$,
%
\begin{eqnarray}
\label{def-tauxi}
\tau^{\varepsilon}_n & \doteq& \inf  \{ t \geq\xi
_{n-1}^\varepsilon \dvtx
d(Z(t), {\mathcaligr V}) = \varepsilon
  \},\nonumber\\[-8pt]\\[-8pt]
\xi^{\varepsilon}_n & \doteq& \inf  \{ t \geq\tau
_n^\varepsilon \dvtx
d(Z(t), {\mathcaligr V}) = \varepsilon/2   \}.\nonumber
\end{eqnarray}
For each $\varepsilon>0$, we have the decomposition
\begin{eqnarray*}
\label{twoterms}
 \sum_{t_i\in\pi^n}|Y(t_i)-Y(t_{i-1})|^p
 & = & \sum
_{t_i\in
\pi^n}\sum_{k=1}^\infty|Y(t_i)-Y(t_{i-1})|^p \mathbb{I}_{
(\tau_k^\varepsilon,\xi_k^\varepsilon)}(t_{i-1}) \nonumber\\
\nonumber
&&{} + \sum_{t_i\in
\pi^n}\sum_{k=0}^\infty|Y(t_i)-Y(t_{i-1})|^p \mathbb{I}_{
[\xi_k^\varepsilon,\tau_{k+1}^\varepsilon]}(t_{i-1}). \nonumber
\end{eqnarray*}
Therefore, for any given $\delta>0$, we have
%
\begin{eqnarray}
\label{Ptwoterms}
&& {\mathbb P}  \biggl(\sum_{t_i\in\pi^n}|Y(t_i)-Y(t_{i-1})|^p
>\delta
\biggr) \nonumber\\[-1pt]
&&\qquad \leq {\mathbb P}  \Biggl(\sum_{t_i\in\pi^n}\sum
_{k=1}^\infty
|Y(t_i)-Y(t_{i-1})|^p \mathbb{I}_{[\tau_k^\varepsilon,\xi
_k^\varepsilon)}(t_{i-1})>\frac
{\delta}{2}  \Biggr) \\[-1pt]
 &&\quad\qquad{} +  {\mathbb P}  \Biggl(\sum_{t_i\in\pi
^n}\sum
_{k=0}^\infty
|Y(t_i)-Y(t_{i-1})|^p \mathbb{I}_{[\xi_k^\varepsilon,\tau
_{k+1}^\varepsilon
)}(t_{i-1})>\frac{\delta}{2}  \Biggr).\nonumber
\end{eqnarray}

Under additional uniform boundedness assumptions on $b$ and $\sigma$,
the proof of~(\ref{YzeroE}) is essentially a consequence of the
following two lemmas, which provide estimates on the two terms on the
right-hand side of (\ref{Ptwoterms}).

\begin{lemma}
\label{lem:Y5}
Suppose $b$ and $\sigma$ are uniformly bounded.
Then, for each $\varepsilon> 0$,
%
\begin{equation}\label{est:Y5}
\lim_{\Delta(\Pi_n) \rightarrow0}
{\mathbb P}  \Biggl(\sum_{t_i\in\pi^n}\sum_{k=1}^\infty
|Y(t_i)-Y(t_{i-1})|^p \mathbb{I}_{
[\tau_k^\varepsilon,\xi_k^\varepsilon)}(t_{i-1})>\frac{\delta
}{2}  \Biggr) = 0.
\end{equation}
\end{lemma}
\begin{pf} Fix $\varepsilon>0$, $n \in{\mathbb N}$, and let
\[
\Omega_n^\varepsilon\doteq
  \biggl\{Z(t) \notin{\mathcaligr V},\  \forall t \in\bigcup
_{k\in{\mathbb N}
:\xi_k^\varepsilon\leq T}   [\xi_k^\varepsilon, \xi
_k^\varepsilon+ \Delta(\pi
^n)  ]
  \biggr\}.
\]
Also, define
\[
N^\varepsilon\doteq\inf  \{k\geq0 \dvtx \mbox{either }\tau
_k^\varepsilon>T
\mbox{ or }\xi_k^\varepsilon>T  \}.
\]
Observe that ${\mathbb P}$-a.s., $N^\varepsilon< \infty$ since
$Z$ has continuous sample paths and therefore crosses the levels $\{z
\in G \dvtx d(z, {\mathcaligr V})= \varepsilon\}$
and $\{z \in G \dvtx d(z, {\mathcaligr V}) =\varepsilon/2\}$ at most a finite
number of
times in the interval $[0,T]$.
The continuity of $Z$ also implies that for each $\varepsilon> 0$,
%
\begin{equation}\label{est:Y6}
{\mathbb P}  (\Omega_n^\varepsilon
  ) \rightarrow1\qquad  \mbox{as }  \Delta(\pi
^n)\rightarrow0.
\end{equation}

On the set $\Omega_n^\varepsilon$, we have
%
\begin{eqnarray}
\label{eq-lemubound}
&& \sum_{t_i\in\pi^n}\sum_{k=1}^\infty
|Y(t_i)-Y(t_{i-1})|^p \mathbb{I}_{
[\tau_k^\varepsilon,\xi_k^\varepsilon)}(t_{i-1}) \nonumber\\[-1pt]
 &&\qquad \leq\max_{t_i \in\pi^n} |Y(t_i)-Y(t_{i-1})|^{p-1}
\sum_{t_i\in\pi^n}\sum_{k=1}^\infty
L(t_{i-1},t_{i}] \mathbb{I}_{[\tau_k^\varepsilon,\xi_k^\varepsilon
)}(t_{i-1})\nonumber\\[-9pt]\\[-9pt]
 &&\qquad = \max_{t_i \in\pi^n} |Y(t_i)-Y(t_{i-1})|^{p-1}
\sum_{t_i\in\pi^n} \sum_{k=1}^\infty L(t_{i-1}, t_{i}] \mathbb{I}_{
[\tau_k^\varepsilon,\xi_k^\varepsilon)}(t_{i-1}) \nonumber\\[-1pt]
 &&\qquad \leq
\max_{t_i \in\pi^n} |Y(t_i)-Y(t_{i-1})|^{p-1}
\sum_{k=1}^\infty L\bigl(\tau_k^\varepsilon\wedge T, \bigl(\xi_k^\varepsilon+
\Delta(\pi^n)\bigr)\wedge T\bigr].\nonumber
\end{eqnarray}
By definition, ${\mathbb P}$-a.s. $(Z,Y)$ satisfy the ESP for $X$.
Therefore, by
Lemma \ref{Lip}, ${\mathbb P}$-a.s., for each $k \in{\mathbb N}$,
$(Z(\tau_k^\varepsilon\wedge T + \cdot), Y(\tau_k^\varepsilon
\wedge T + \cdot) -
Y(\tau_k^\varepsilon\wedge T))$
solve the ESP for
$Z(\tau_k^\varepsilon\wedge T) + X(\tau_k^\varepsilon\wedge T +
\cdot) - X(\tau
_k^\varepsilon\wedge T)$.
On $\Omega^{\varepsilon}_n$, $Z$ is away from ${\mathcaligr V}$ on
$[\tau
_k^\varepsilon
\wedge T, (\xi_k^\varepsilon+ \Delta(\pi^n))\wedge T]$
for each $k\geq1$, and hence by Theorem 2.9 of \cite{ram1} it follows that
$L(\tau_k^\varepsilon\wedge T, (\xi_k^\varepsilon+ \Delta(\pi
^n))\wedge T] <
\infty$.
Together with the fact that ${\mathbb P}$-a.s. $N^\varepsilon<\infty
$, this
implies that
\[
\sum_{k=1}^\infty L\bigl(\tau_k^\varepsilon\wedge T, \bigl(\xi_k^\varepsilon+
\Delta(\pi^n)\bigr)\wedge T\bigr] <\infty \qquad\mbox{ ${\mathbb P}$-a.s. on }\Omega
^\varepsilon_n.
\]
On the other hand, since $Y$ is continuous on $[0,T]$ and $p > 1$, we have
\[
\max_{t_i \in\pi^n} |Y(t_i)-Y(t_{i-1})|^{p-1} \rightarrow0
\qquad\mbox{as }\Delta(\pi^n)\rightarrow0.
\]
Combining the above two displays
with (\ref{est:Y6}), we conclude that for every $\delta> 0$, as
$\Delta(\pi^n)\rightarrow0$,
\[
{\mathbb P}  \Biggl( \max_{t_i \in\pi^n}|Y(t_i)-Y(t_{i-1})|^{p-1} \sum
_{k=1}^\infty L\bigl(\tau_k^\varepsilon\wedge T, \bigl(\xi_k^\varepsilon+
\Delta(\pi^n)\bigr)\wedge T\bigr] > \frac{\delta}{2}   \Biggr) \rightarrow0.
\]
Together with (\ref{eq-lemubound}), this shows that (\ref{est:Y5})
holds and completes the proof of the lemma.
\end{pf}

In the next lemma, $q \geq2$ is the value for which
Assumption \ref{as-weakcont} is satisfied.

\begin{lemma}
\label{lem:Y4}
Suppose $b$ and $\sigma$ are uniformly bounded.
Then there exists a finite constant $C < \infty$ such that for each
$\varepsilon> 0$,
%
\begin{eqnarray}
\label{est:Y4-pvar}
 &&\lim_{\triangle(\pi^n)\rightarrow0}{\mathbb P}
\Biggl(\sum
_{t_i\in\pi
^n}\sum_{k=0}^\infty
|Y(t_i)-Y(t_{i-1})|^p\mathbb{I}_{[\xi_k^\varepsilon,\tau
_{k+1}^\varepsilon
)}(t_{i-1})>\frac{\delta}{2}   \Biggr)\nonumber\\[-8pt]\\[-8pt]
   &&\qquad
   \leq
\cases{\displaystyle
 \frac{C}{\delta}{\mathbb E}  \Biggl[\int_0^{T} \sum_{k=0}^\infty\mathbb{I}_{[\xi_k^\varepsilon,\tau
_{k+1}^\varepsilon]}(t) \,dt  \Biggr],
&\quad  if  $q = 2$, \vspace*{2pt}\cr
0, &\quad  if $q >2$.
}\nonumber
\end{eqnarray}
\end{lemma}
\begin{pf}
Fix $\varepsilon> 0$. Then
by Markov's inequality [whose application is justified by (\ref{finite-ineq})
when $q > 2$, and by (\ref{just-Markov}) when $q=2$]
and the monotone convergence theorem
%
\begin{eqnarray}
\label{eq-Markov}
&& {\mathbb P}  \Biggl(\sum_{t_i\in\pi^n}\sum_{k=0}^\infty
|Y(t_i)-Y(t_{i-1})|^p\mathbb{I}_{[\xi_k^\varepsilon,\tau
_{k+1}^\varepsilon
)}(t_{i-1})>\frac{\delta}{2}   \Biggr) \nonumber\\[-8pt]\\[-8pt]
 &&\qquad \leq\frac{2}{\delta} \sum_{t_i\in\pi^n}\sum
_{k=0}^\infty{\mathbb E} \bigl [
|Y(t_i)-Y(t_{i-1})|^p\mathbb{I}_{[\xi_k^\varepsilon,\tau
_{k+1}^\varepsilon)}(t_{i-1})
  \bigr].\nonumber
\end{eqnarray}
Recall that $a = \sigma^T \sigma$, and let $\bar{C} > 1$ be an upper
bound on $|b|$,
$|\sigma|$ and $|a|$.
By Assumption \ref{as-weakcont}, the definition (\ref{def-x}) of $X$
and the elementary inequality
$|x+y|^q \leq2^q(|x|^q + |y|^q)$, there exists $K_T < \infty$ such
that for each $t_i \in\pi^n$,
%
\begin{eqnarray*}
 &&{\mathbb E}  [ |Y(t_i)-Y(t_{i-1})|^p | {\mathcaligr F}_{t_{i-1}}  ] \\
&& \qquad\leq  K_T {\mathbb E}  \Bigl[ \sup_{u \in[t_{i-1},t_i]}
  |X (u)
- X(t_{i-1})   |^q \big| {\mathcaligr F}_{t_{i-1}}   \Bigr] \\
&& \qquad  \leq2^q K_T {\mathbb E} \biggl [\sup_{u \in[t_{i-1},t_i]}
\bigg|\int_{t_{i-1}}^u b(Z(v)) \, dv   \bigg|^q   \\
&& \qquad \hspace*{47pt}{}+
\sup_{u \in[t_{i-1}, t_i]}  \bigg |\int_{t_{i-1}}^{u} \sigma(Z(v))
\, dB_v  \bigg |^q  \Big| {\mathcaligr F}_{t_{i-1}}  \biggr]
\\
 && \qquad \leq2^q K_T {\mathbb E}  \biggl[ \bar{C}^q(t_{i} -
t_{i-1})^q +
  \biggl(\frac{q}{q-1}  \biggr)^q  \bigg |\int_{t_{i-1}}^{t_i} \sigma
(Z(v)) \, dB_v  \bigg |^q \Big | {\mathcaligr F}_{t_{i-1}}  \biggr]
\\
&& \qquad  \leq2^q K_T \bar{C}^q (t_{i} - t_{i-1})^q \\
&& \qquad\quad  {}+
2^q K_T  \biggl(\frac{q}{q-1}  \biggr)^q \tilde K{\mathbb E}  \biggl[
\biggl(\int
_{t_{i-1}}^{t_i}   | a(Z(v))   | \, dv  \biggr)^{q/2}  \Big|
{\mathcaligr F}_{t_{i-1}}
  \biggr],
\end{eqnarray*}
where the third inequality holds due to the uniform bound on $b(\cdot
)$, the
Markov property of $Z$ and Doob's maximal martingale inequality,
while the fourth inequality follows, with $\tilde{K} < \infty$ a universal
constant, by an application of the martingale moment inequality,
which is justified since the uniform boundedness on $a$ ensures that
the stochastic integral is a \vspace*{1pt}martingale.

Define $\tilde C \doteq2^q K_T [\bar{C}^q \vee(q^q \bar{C}^{q/2}
\tilde{K}/(q-1)^q)]$. Using the bound on $a$, the last inequality
shows that for each $t_i \in\pi^n$,
%
\begin{equation}
\label{new-ineqset}
{\mathbb E}  [ |Y(t_i)-Y(t_{i-1})|^p  | {\mathcaligr
F}_{t_{i-1}}  ]
\leq\tilde{C}   [ (t_i - t_{i-1})^q + (t_i - t_{i-1})^{q/2}
  ].
\end{equation}
We now consider two cases. If $q>2$, it follows from (\ref
{new-ineqset}) that,
for all sufficiently large $n$ such that $\Delta(\pi^n) < 1$,
\[
{\mathbb E}  [ |Y(t_i)-Y(t_{i-1})|^p \vert {\mathcaligr F}_{t_{i-1}}
] \leq
2 \tilde
C \Delta(\pi^n)^{q/2-1} (t_{i} - t_{i-1}).
\]
Multiplying both sides of this inequality by
$\mathbb{I}_{[\xi_k^\varepsilon,\tau_{k+1}^\varepsilon)} (t_{i-1})$,
which is ${\mathcaligr F}_{t_{i-1}}$-measurable since
$\tau_k^\varepsilon$ and $\xi_k^\varepsilon$ are stopping times, then
taking expectations and subsequently
summing over $k=0, 1, \ldots ,$ and $t_i \in\pi^n$,
it follows that
%
\begin{equation}
\label{finite-ineq}
 \sum_{t_i\in\pi^n}\sum_{k=0}^\infty{\mathbb E}  \bigl[
|Y(t_i)-Y(t_{i-1})|^p \mathbb{I}_{[\xi_k^\varepsilon,\tau
_{k+1}^\varepsilon)}(t_{i-1})
  \bigr]
\leq
2\tilde{C} \Delta(\pi^n)^{q/2-1} T.
\end{equation}
Since $\Delta(\pi^n)^{q/2-1} \rightarrow0$ as
$n \rightarrow\infty$,
combining this with (\ref{eq-Markov}), we then obtain
\[
\lim_{\triangle(\pi^n)\rightarrow0}{\mathbb P}  \Biggl(\sum_{t_i\in
\pi
^n}\sum_{k=0}^\infty
|Y(t_i)-Y(t_{i-1})|^p\mathbb{I}_{[\xi_k^\varepsilon,\tau
_{k+1}^\varepsilon
)}(t_{i-1})>\frac{\delta}{2}
  \Biggr) =0.
\]

On the other hand, if $q=2$, again
multiplying both sides of (\ref{new-ineqset}) by $\mathbb{I}_{[\xi
_k^\varepsilon
,\tau_{k+1}^\varepsilon)} (t_{i-1})$, then
taking expectations, subsequently summing over $k=0, 1, \ldots ,$ and
$t_i \in
\pi^n$, and then using the monotone
convergence theorem to interchange expectation and summation, we obtain
%
\begin{eqnarray*}
\label{just-Markov}
&& \sum_{t_i\in\pi^n}\sum_{k=0}^\infty{\mathbb E} \bigl [
|Y(t_i)-Y(t_{i-1})|^p \mathbb{I}_{[\xi_k^\varepsilon,\tau
_{k+1}^\varepsilon)}(t_{i-1})
  \bigr] \\
 && \qquad
 \leq \tilde{C}   \Biggl( \Delta(\pi^n)^{q} +
{\mathbb E} \Biggl [\sum_{t_i\in\pi^n}(t_{i} - t_{i-1}) \sum
_{k=0}^\infty
\mathbb{I}_{[\xi_k^\varepsilon, \tau_{k+1}^\varepsilon)} (t_{i-1})
  \Biggr]   \Biggr) <
\infty .
\end{eqnarray*}
%
Sending $\Delta(\pi^n) \rightarrow0$ on both sides of this inequality
and invoking
the bounded convergence theorem, the right-continuity of $\mathbb
{I}_{[\xi
_k^\varepsilon, \tau_{k+1}^\varepsilon)} (\cdot)$ and the definition
of the Riemann integral, we obtain
\begin{eqnarray*}
&&\lim_{\Delta(\pi^n) \rightarrow0} \sum_{t_i \in\pi^n}\sum
_{k=0}^\infty
{\mathbb E}  \bigl[
|Y(t_i)-Y(t_{i-1})|^p \mathbb{I}_{[\xi_k^\varepsilon,\tau
_{k+1}^\varepsilon)}(t_{i-1})
  \bigr] \\
&& \qquad \leq\tilde{C} {\mathbb E}  \Biggl[ \int_0^{T}
\sum_{k=0}^\infty\mathbb{I}_{[\xi_k^\varepsilon,\tau
_{k+1}^\varepsilon)}(t) \, dt
  \Biggr].
\end{eqnarray*}
Together with (\ref{eq-Markov}), this shows that (\ref{est:Y4-pvar})
holds with $C = 2 \tilde{C}$.
\end{pf}

\begin{pf*}{Proof of Theorem \ref{th-main2}}
Due to Lemma \ref{lem:Ytau}, using a localization argument and the
local boundedness
of $b$ and $\sigma$ stated in Assumption \ref{as-locbd},
we can assume without loss of generality that $a, b$ and $\sigma$ are bounded.
Then, combining (\ref{Ptwoterms}) with Lemmas \ref{lem:Y5} and \ref
{lem:Y4}, we have
\begin{eqnarray*}
 &&\lim_{\Delta(\pi^n)\rightarrow0}{\mathbb P}
\biggl(\sum_{t_i\in
\pi^n}|Y(t_i)-Y(t_{i-1})|^p>\delta  \biggr) \\
 &&\qquad \leq
\cases{
 \displaystyle\frac{C}{\delta}
{\mathbb E}  \Biggl[\int_0^{T} \sum_{k=0}^\infty\mathbb{I}_{[\xi
_k^\varepsilon
,\tau
_{k+1}^\varepsilon]}(t) \,dt  \Biggr],
&\quad   if   $q = 2$, \vspace*{2pt}\cr
0, &\quad  if $q >2$,
}
\end{eqnarray*}
for every $\varepsilon> 0$, and so (\ref{YzeroE}) holds for the case
$q > 2$.
If $q = 2$, sending $\varepsilon\downarrow0$ and using the bounded
convergence theorem and the definition
of the stopping times $\xi_k^\varepsilon$ and $\tau_k^\varepsilon$,
we see that
the term on the right-hand side converges to
\[
\frac{C}{\delta} {\mathbb E}  \biggl[ \int_0^{T} \mathbb
{I}_{{\mathcaligr V}}
(Z(t)) \,
dt   \biggr] = 0,
\]
where the last equality follows from (\ref{eq-zero}) and the fact that
${\mathcaligr V} \subset\partial G$.
This proves (\ref{YzeroE}), and
Theorem \ref{th-main2} then follows from the
discussion at the beginning of Section \ref{sec-dirichlet}.
\end{pf*}

\subsection{\texorpdfstring{Proof of Corollary
\protect\ref{cor-egburtob}}{Proof of Corollary 3.7}}
\label{subs-pfegburtob}

Suppose that, as in (\ref{lrcond}),
the functions $L$ and $R$ on $[0,\infty)$ are given by
$L(y) = -c_L y^{\alpha_L}$ and $R(y) = c_R y^{\alpha_R}$ for some
$\alpha_L, \alpha_R, c_L, c_R \in(0, \infty)$.
As defined in Section \ref{sec-egburtob},
let $(G,d(\cdot))$ and $\overline{\Gamma}$ be the associated
ESP and ESM, and let $Z = (Z_1, Z_2)$
be the associated two-dimensional RBM: $Z = \overline{\Gamma}(B)$,
where $B = (B_1, B_2)$ is a standard two-dimensional Brownian motion.
Then Assumptions \ref{as-weak}
and \ref{as-locbd} are automatically satisfied for this
family of reflected diffusions.
In order to prove the corollary, it suffices to show that
Assumption \ref{as-weakcont} holds.
Indeed, then all the assumptions of Theorem \ref{th-main2} are
satisfied, and Corollary \ref{cor-egburtob} follows as a
consequence.

We now recall the representation for $Y \doteq Z - B$ that
was obtained in Section 4.3 of \cite{burkanram08}.
First, note that $Z_2$ is a one-dimensional RBM on
$[0,\infty)$ with the pathwise representation $Z_2 = \Gamma_1(B_2)$,
where
$\Gamma_1$ is the one-dimensional reflection map on $[0,\infty)$.
Thus, $Y_2 = \Lambda_2 (B_2)$, where $\Lambda_2 (\psi) \doteq\Gamma
_1(\psi) - \psi$
is
given explicitly by
%
\begin{equation}
\label{def-Lambda}
\Lambda_2(\psi)(t) \doteq\sup_{0\leq s\leq t}[-\psi(s)]^+,\qquad
\psi\in
{\mathcaligr C}[0,\infty), t \in[0,\infty).
\end{equation}
(Recall that ${\mathcaligr C}[0,\infty)$ is
the space of continuous functions on $[0,\infty)$,
equipped with the topology of uniform convergence on compact
sets.)
Since $Y_2$ is a nondecreasing process, it is clearly of finite
variation. Therefore, to establish Assumption \ref{as-weakcont},
it suffices to show that the inequality (\ref{l2cont}) holds
with $Y$ replaced by $Y_1$.
From Section 4.3 of
\cite{burkanram08}, it follows that pathwise
$Z_1 = \bar{\Gamma}_{\ell,r} (B_1)$,
where $\bar{\Gamma}_{\ell,r}$ is the ESM whose domain is the time-dependent
interval $[l(\cdot), r(\cdot)]$, with $l(t) \doteq L(Z_2(t))$ and
$r(t)\doteq R(Z_2(t))$, for $t \in[0,\infty)$. A precise definition of
$\bar{\Gamma}_{\ell,r}$ is stated as Definition 2.2 of \cite{burkanram08},
but for the present purpose it suffices to note that Theorem 2.6 of
\cite{burkanram08} establishes
the explicit representation $\bar{\Gamma}_{\ell,r} (\psi) = \psi-
\Xi_{\ell,r}(\psi)$,
where for $\psi\in{\mathcaligr C}[0,\infty)$
such that $\psi(0) \in[\ell(0), r(0)]$, and
$t \in[0,\infty)$,
%
\begin{eqnarray}
\label{def-Xi}
&&\Xi_{\ell,r}(\psi)(t) \doteq \max
  \Bigl(  \Bigl [0 \wedge\inf_{u\in
[0,t]}\bigl(\psi(u)-\ell(u)\bigr)
  \Bigr],   \nonumber\\[-8pt]\\[-8pt]
&&\hspace*{86.3pt}\sup_{s\in[0,t]}  \Bigl[\bigl(\psi(s)-r(s)\bigr)\wedge\inf_{u\in
[s,t]}\bigl(\psi(u)-\ell(u)\bigr)  \Bigr]   \Bigr).\nonumber
\end{eqnarray}
Thus, we see that $Y_1 = \Lambda_1(B_1,B_2)$, where
$\Lambda_1$ is the map from
${\mathcaligr C}[0,\infty)^2$ to ${\mathcaligr C}[0,\infty)$ given by
\[
\Lambda_1 \dvtx (\psi_1, \psi_2) \mapsto
-\Xi_{L \circ\Gamma_1 (\psi_2), R\circ\Gamma_1 (\psi_2)} (\psi
_1).
\]
{F}rom the explicit expression for $\Xi_{\ell,r}$ given in (\ref{def-Xi}),
it can be easily verified
that the map $(\ell,r,\psi) \mapsto-\Xi_{\ell,r}(\psi)$ from
${\mathcaligr C}[0,\infty)^3$ to ${\mathcaligr C}[0,\infty)$ is
Lipschitz continuous.
In addition, it follows from (\ref{def-Lambda}) that the map $\psi
\mapsto
\Gamma_1 (\psi)$ from ${\mathcaligr C}[0,\infty)$ to itself is also
Lipschitz continuous.
If $L$ and $R$ are H\"{o}lder continuous with
exponent $\alpha= \alpha_L \wedge\alpha_R \in(0,1)$, it follows
that the
composition maps $\ell= L \circ\Gamma_1$ and $r = R \circ\Gamma_1$
are also H\"{o}lder continuous with exponent $\alpha$.
When combined, the above statements then imply that the map $\Lambda
_1$ is locally
H\"{o}lder continuous on ${\mathcaligr C}[0,\infty)^2$ with exponent
$\alpha
$, and so
(\ref{l2cont}) holds for any
$p \geq2/\alpha$ with, correspondingly, $q = \alpha p$.
On the other hand, if $L$ and
$R$ are locally Lipschitz continuous [i.e., if (\ref{lrcond})
is satisfied with
$\alpha\geq1$], then clearly $\bar{\Gamma}$ is also locally Lipschitz
continuous, and
so~(\ref{l2cont}) holds with $p = q = 2$. Thus, the result
follows in this case as well (note that, due to the localization
result of Section \ref{subs-equiv}, it suffices for the ESM to be locally
Lipschitz or locally H\"{o}lder, that is, Lipschitz continuous or H\"
{o}lder continuous on paths that lie in a compact
set on any finite time interval).

\begin{remark}
\label{rem54}
The proof above also shows that when $\alpha < 1$, $Z =
B + A$, where $A$ is a process of zero $p$-variation for every $p > 2/\alpha$.
However, this is likely to be a sub-optimal result, since given that
$Z$ is a Dirichlet process even when the domain is cusp-like, one would
expect that $Z$ would also be a Dirichlet process when the domain is
flatter (corresponding to $\alpha > 1$). Indeed, more generally, it would
be of interest to determine the lowest $p$-variation that vanishes for a
given~$\alpha$, to better understand the relationship between the ``roughness''
of the paths of $Z$ and the curvature of the boundary of the
domain. Such questions motivate a ``rough paths'' analysis (see, e.g.,
\cite{LyonsStFlour04} and \cite{FrizVictBook09}) of reflected stochastic processes.
\end{remark}

\begin{remark}
\label{rem-egburtob}
The above class of reflected diffusions provides one example
of a situation where the ESM is locally H\"{o}lder continuous, but the
(generalized) completely-${\mathcaligr S}$ condition does not hold.
However, we believe that it should be possible to combine a
localization argument of the kind
used in \cite{dupish93} and the sufficient condition for Lipschitz
continuity of the ESM obtained in
Theorem 3.3 of \cite{ram1} to identify a broad class of
piecewise smooth domains and directions of reflection where the
generalized completely-${\mathcaligr S}$ condition fails to hold, but
for which the associated ESM is locally H\"{o}lder continuous.
\end{remark}

\begin{appendix}\label{appm}
\section*{Appendix A: Elementary properties of the ESP}
\setcounter{theorem}{0}
\begin{lemma}
\label{Lip}
If $(\phi,\eta)$ is a solution to the ESP $(G,d(\cdot))$ for $\psi
\in{\mathcaligr C}_G [0,\infty)$, then for
each $0 \leq s < \infty$, $(\phi^s,\eta^s)$ is a
solution to the ESP for $\phi(s)+\psi^s$, where $\phi^s(\cdot
)\doteq\phi(s+\cdot)$,
\[
\psi^s(\cdot)\doteq\psi(s+\cdot)-\psi(s)  \quad\mbox{and}\quad
 \eta^s(\cdot)\doteq\eta(s+\cdot)-\eta(s).
\]
Moreover, if the ESM is well-defined and Lipschitz continuous on
${\mathcaligr C}_G [0,\infty)$ then for every $T < \infty$,
there exists $\tilde{K}_T < \infty$ such that for every
$0\leq s<t \leq T+s$,
\[
|\eta(t)-\eta(s)| \leq\tilde{K}_T \sup_{u \in[0, t-s]}|\psi
(s+u)-\psi(s)|.
\]
\end{lemma}
\begin{pf} Fix $s \in[0,\infty)$ and a path $\psi\in{\mathcaligr
D}_G [0,\infty)$.
The first statement follows from Lemma 2.3 of \cite{ram1}. It implies that
$\eta^s=\bar{\Gamma}(\psi^1)-\psi^1$, where $\psi^1 \doteq\phi(s)
+\psi^s$.
On the other hand, consider the path $\psi^2$ which is
equal to the constant $\phi(s)$ on $[0,\infty)$, that is, $\psi^2 (u)
\doteq\phi(s)$ for all $u \in[0,\infty)$.
Then clearly $(\phi(s),0)$ is the
unique solution to the ESP for $\psi^2$, that is, $0 = \bar{\Gamma
}(\psi^2)
(u) - \psi^2(u)$ for all $u \in[0,\infty)$.
Using the Lipschitz continuity of the ESM, for $\delta\in[0,T-s]$ we obtain
%
\begin{eqnarray*}
|\eta^s(\delta)-0|&\leq& \sup_{u\in[0, \delta]}  |\bar
{\Gamma}
(\psi^1)(u)-\psi^1 (u) -
\bar{\Gamma}(\psi^2)(u)+\psi^2(u)  | \\
&\leq& \sup_{u\in[0, \delta]}
  |\bar{\Gamma}(\psi^1)(u)- \bar{\Gamma}(\psi^2)(u)  | +
\sup_{u\in[0, \delta]} |\psi^1(u) - \psi^2(u)|\\
&\leq& K_T \sup_{u\in[0, \delta]}  |\psi^s(u)  |+ \sup
_{u\in[0, \delta]} |\psi^s(u)|,
\end{eqnarray*}
where $K_T < \infty$ is the Lipschitz constant of $\bar{\Gamma}$ on $[0,T]$.
The lemma follows by letting $\tilde{K}_T \doteq K_T +1$ and $\delta=
t-s$.
\end{pf}

\section*{Appendix B: Auxiliary results}

For completeness, we provide the proof of the fact that the sequences
of times
defined in Section \ref{subsub-auxres} are stopping times.
\renewcommand{\thetheorem}{B.\arabic{theorem}}
\setcounter{theorem}{0}
\renewcommand{\thelemma}{B.\arabic{lemma}}

\begin{lemma}
\label{lem-stop}
$\{\newstop^\varepsilon_n\}_{n \in{\mathbb N}}$, $\{\newstop
_{(k),n}^\varepsilon\}
_{n \in
{\mathbb N}}$, $k \in{\mathbb N}$, are sequences of
$\{{\mathcaligr F}_t\}$-stopping times. Also, $\{\newstop
^{k,\varepsilon}_n\}
_{n \in
{\mathbb N}}$, $k \in{\mathbb N}$, are sequences of
$\{{\mathcaligr F}^k_t\}$-stopping times.
\end{lemma}
\begin{pf} Clearly,
$\newstop_0^\varepsilon\doteq0$ is an $\{{\mathcaligr F}_t\}
$-stopping time.
Now, suppose $\newstop_{n-1}^\varepsilon$ is an $\{{\mathcaligr F}_t\}
$-stopping
time and
note that for each $\varepsilon> 0$, $n \in{\mathbb N}$ and $t \in
[0,\infty)$,
\[
\{\newstop_n^\varepsilon\leq t \} = \bigcup_{k \in{\mathbb Z}}
  [
\{\newstop_{n-1}^\varepsilon\leq t\} \cap  \{ Z(\newstop
_{n-1}^\varepsilon)
\in H_{2^k \varepsilon}   \} \cap A_{k,n}^\varepsilon(t)  ],
\]
where
\[
A_{k,n}^\varepsilon(t) \doteq
 \Bigl \{ \sup_{s\in[\newstop_{n-1}^\varepsilon, t]} \langle Z(s),
\vec{{\mathbf v}}
\rangle\geq2^{k+1} \varepsilon  \Bigr\}
\cup  \Bigl\{ \inf_{s \in[\newstop_{n-1}^\varepsilon, t]} \langle Z(s),
\vec{{\mathbf v}}\rangle\leq2^{k-1}\varepsilon  \Bigr\}.
\]
Then $\{\newstop_{n-1}^\varepsilon\leq t \} \in{\mathcaligr F}_t$ because
$\newstop
_{n-1}^\varepsilon$ is an $\{{\mathcaligr F}_t\}$-stopping
time. Since $Z$ is continuous
we also know that $\{\newstop_{n-1}^\varepsilon\leq t\} \cap  \{
Z(\newstop_{n-1}^\varepsilon) \in H_{2^k \varepsilon}   \}$
lies in ${\mathcaligr F}_t$. In addition, the continuity of $\langle Z,
\vec{{\mathbf v}}
\rangle$ and the fact that
$[2^{k+1} \varepsilon, \infty)$ and $(-\infty, 2^{k-1} \varepsilon
]$ are closed
show that
$\{\newstop_{n-1}^\varepsilon\leq t \} \cap A_{n,k}^\varepsilon
(t)\in{\mathcaligr F}_t$.
When combined, this implies that $\{\newstop_n^\varepsilon\leq t\}
\in
{\mathcaligr
F}_t$ or,
equivalently, that $\newstop_n^\varepsilon$ is an $\{{\mathcaligr
F}_t\}
$-stopping time,
and the first assertion follows by induction. The proof for the other
sequences is exactly analogous.
\end{pf}
\end{appendix}

\section*{Acknowledgments}

The authors are grateful to F. Coquet for first posing the question
as to whether the generalized processor sharing reflected diffusion is
a Dirichlet process.
The second author would also like to thank S. R. S. Varadhan for his hospitality
and advice during
her stay at the Courant Institute, during which part of this work was
completed, and is
grateful to A. S.-Sznitman for a useful discussion.


\printaddresses

\end{document}